\documentstyle{amltd2004}
\begin{document}
\annalsline{156}{2002}
\received{March 30, 2000}
\startingpage{797}
\def\bye{\end{document}}
 \font\tenrm=cmr10
\input amssym.def
\input amssym.tex
\def\ritem#1{\item[{\rm #1}]}

%--------------- Author macros ---------------
%for Bbb in amstex
\catcode`\@=11
\font\twelvemsb=msbm10 scaled 1100
\font\tenmsb=msbm10
%\font\ninemsb=msbm7 scaled 1100%msbm9
\font\ninemsb=msbm10 scaled 800
\newfam\msbfam
\textfont\msbfam=\twelvemsb  \scriptfont\msbfam=\ninemsb
  \scriptscriptfont\msbfam=\ninemsb
\def\msb@{\hexnumber@\msbfam}
\def\Bbb{\relax\ifmmode\let\next\Bbb@\else
 \def\next{\errmessage{Use \string\Bbb\space only in math
mode}}\fi\next}
\def\Bbb@#1{{\Bbb@@{#1}}}
\def\Bbb@@#1{\fam\msbfam#1}
\catcode`\@=12

 \catcode`\@=11
\font\twelveeuf=eufm10 scaled 1100
\font\teneuf=eufm10
\font\nineeuf=eufm7 scaled 1100%eufm9
\newfam\euffam
\textfont\euffam=\twelveeuf  \scriptfont\euffam=\teneuf
  \scriptscriptfont\euffam=\nineeuf
\def\euf@{\hexnumber@\euffam}
\def\frak{\relax\ifmmode\let\next\frak@\else
 \def\next{\errmessage{Use \string\frak\space only in math
mode}}\fi\next}
\def\frak@#1{{\frak@@{#1}}}
\def\frak@@#1{\fam\euffam#1}
\catcode`\@=12
%-------------- Author entries --------------------

\title{An infinite Ramsey theorem and\\ some Banach-space dichotomies} %Article title
\shorttitle{An infinite Ramsey theorem}   
 \author{W. T. Gowers}
 \institutions{DPMMS, University of Cambridge, Wilberforce Road, Cambridge CB3 OWB, UK\\
{\eightpoint {\it E-mail address\/}: wtg10@dpmms.cam.ac.uk}}
%-------------- Article Text--------------------

%\intro %(Optional, Introduction)

\font\deffont=cmssi10
\def\sqr{$\vcenter{\hrule height .3mm
\hbox {\vrule width .3mm height 2mm \kern 2mm
\vrule width .3mm} \hrule height .3mm}$}

\def \casespace{\noalign{\vskip 5 pt}}
\def \R{{\Bbb R}}
\def \E{{\Bbb E}}
\def \P{{\Bbb P}}
\def \M{{\Bbb M}}
\def \N{{\Bbb N}}
\def \C{{\Bbb C}}
\def \D{\Delta}
\def \normof#1{\left\|#1\right\|}
\def \nm#1{\left\|#1\right\|}
\def \bnm#1{\bigl\|#1\bigr\|}
\def \Bnm#1{\Bigl\|#1\Bigr\|}
\def \bgnm#1{\biggl\|#1\biggr\|}
\def \nmm#1{\left|\!\left|\!\left|#1\right|\!\right|\!\right|}
\def \bnmm#1{\bigl|\!\bigl|\!\bigl|#1\bigr|\!\bigr|\!\bigr|}
\def \Bnmm#1{\Bigl|\!\Bigl|\!\Bigl|#1\Bigr|\!\Bigr|\!\Bigr|}
\def \bgnmm#1{\biggl|\!\biggl|\!\biggl|#1\biggr|\!\biggr|\!\biggr|}
\def \seq#1#2{#1_1,\dots,#1_#2}
\def \sleq#1#2{#1_1<\dots<#1_#2}
\def \speq#1#2{#1_1+\dots+#1_#2}
\def \sm#1#2{\sum_{#1=1}^#2}
\def \ge{\geqslant}
\def \le{\leqslant}
\def \e{\varepsilon}
\def \G{\Gamma}
\def \g{\gamma}
\def \d{\delta}
\def \t{\tau}
\def \m{\mu}
\def \n{\nu}
\def \th{\theta}
\def \supp{\mathop{\rm supp}}
\def \sign{{\rm sign}}
\def \a{\alpha}
\def \b{\beta}
\def \l{\lambda}
\def \c{\choose}
\def \o{\over}
\def \sp#1{\langle#1\rangle}
\def \ra{\rightarrow}
\def \om{{\omega}}
\def \s{\sigma}

 \centerline{\bf Abstract}
\vglue12pt

A problem of Banach asks whether every 
infinite-dimensional Banach space which is isomorphic to all its 
infinite-dimensional subspaces must be isomorphic to a separable
Hilbert space. In this paper we prove a result of a Ramsey-theoretic
nature which implies an interesting dichotomy for subspaces of
Banach spaces. Combined with a result of Komorowski and 
Tomczak-Jaegermann, this gives a positive answer to Banach's
problem. We then generalize the Ramsey-theoretic result and
deduce a further dichotomy for Banach spaces with an unconditional
basis.

\section{Introduction}

This paper contains a complete proof of a result announced in [G1]
which has been circulating in preprint form for several years [G2,3].
Our main theorem, when combined with a very different result of
Komorowski and Tomczak-Jaegermann, solves a problem from the famous
1932 book of Banach [B]. He asked whether a separable Hilbert space is
the only infinite-dimensional Banach space, up to isomorphism, which
is isomorphic to every infinite-dimensional closed subspace of
itself. The answer turns out to be yes, but most of the proof appears
to have little to do with the problem. Therefore, in order to motivate
the rest of the paper, we shall begin by explaining (in this section
and the next) how Banach's problem can be reduced to a question with a
much more combinatorial flavour.

To save writing, let us assume from now on that all Banach spaces and
subspaces are infinite-dimensional unless they are specified
otherwise. Recall that a (\/{\it Schauder\/{\rm )} basis} of a Banach space $X$
is a sequence $(x_n)_{n=1}^\infty$ such that every $a\in X$ can be
written uniquely as a norm-convergent sum $a=\sm n \infty a_nx_n$ for
some sequence $(a_n)_{n=1}^\infty$ of scalars. If every $x_n$ has norm
one then the basis is {\it normalized}. It can be shown that if
$(x_n)_{n=1}^\infty$ is a basis, then there is a constant $C$ such
that all the projections $P_N:\sm n \infty a_nx_n \mapsto \sm n N
a_nx_n$ have norm at most $C$. The smallest constant with this
property is called the {\it basis constant} of the basis.  If the
basis constant is 1, then the basis is called {\it monotone}. A basis
$(x_n)_{n=1}^\infty$ is {\it unconditional} if the sum $\sm n \infty
a_nx_n$ converges unconditionally whenever it converges. This, it
turns out, implies that there is a constant $C$ such that, for every
subset $A\subset\N$ (it is sufficient to consider only finite subsets),
the projection $P_A:\sm n \infty a_nx_n\mapsto\sum_{n\in A}a_nx_n$ has
norm at most $C$. The smallest constant with this property is called
the {\it unconditional constant} of the basis.  A sequence
$(x_n)_{n=1}^\infty$ which is a basis/unconditional basis of the
closed subspace that it spans is called a basic sequence/unconditional
basic sequence.

The result of Komorowski and Tomczak-Jaegermann we shall use is the
following [K~T-J].

\proclaim{Theorem} \hskip-9pt Let $X$ be a Banach space with cotype $q$ for
some \hbox{$q<\infty$.} Then either $X$ has a subspace without an unconditional
basis or $X$ has a subspace isomorphic to $\ell_2$.
\endproclaim

Let us follow usual practice and define a Banach space to be {\it
homogeneous} if it is isomorphic to all its subspaces. It can be shown
that every homogeneous Banach space satisfies the cotype condition
above. Indeed, an old result of Mazur states that every Banach space
has a subspace with a basis. On the other hand, Szankowski has shown [Sz]
(generalizing Enflo's solution of the basis problem) that if $p<2$ and
$q>2$, then any Banach space which either fails to have type $p$ or
fails to have cotype $q$ has a subspace without a basis.  Hence, such
a space cannot be homogeneous. So an immediate consequence of Theorem
1 and the definition of homogeneity is the following additional result
of Komorowski and Tomczak-Jaegermann.

\proclaim{{C}orollary} Let $X$ be a homogeneous Banach space. Then 
either $X$ is isomorphic to $\ell_2$ or $X$ fails to have an
unconditional basis. 
\endproclaim

Notice that in the second case no subspace of $X$ has an unconditional
basis, since $X$ is homogeneous. Thus Corollary 1.2 appears to imply a
very strong property of homogeneous spaces not isomorphic to
$\ell_2$. Indeed, it is not at all obvious that there exists a Banach
space such that no subspace has an unconditional basis. 

As it happens, the existence of such spaces was itself a long-standing
open problem until 1991, when counterexamples were found [GM1] (see
also [AD], [F1], [G6], [GM2], [H]). However, these counterexamples did not come close to
being homogeneous. In fact, many of them had a property which is
almost the reverse of homogeneity. A Banach space $X$ is called {\it
decomposable} if it can be written as a direct sum $Y+Z$ with the
projections to $Y$ and $Z$ continuous. Equivalently, $X$ is
decomposable if it admits a nontrivial projection (that is, one of
infinite rank and infinite corank). Otherwise, $X$ is {\it
indecomposable}. It is {\it hereditarily indecomposable} if it has no
decomposable subspace. This is the remarkable property which was
enjoyed by many of the counterexamples. Later we shall have more to
say about hereditarily indecomposable spaces, but for now we shall
simply quote the following result from [GM1].

\proclaim{Theorem} A hereditarily indecomposable Banach space
is isomorphic to no proper subspace of itself.
\endproclaim

 In particular, a hereditarily indecomposable space is 
certainly not homogeneous.

We can now state our first Banach-space dichotomy, which is one of the
main results of this paper.

\proclaim{Theorem}  Every Banach space $X$ has a subspace $W$
which either has an unconditional basis or is hereditarily indecomposable.
\endproclaim

Note that this is a genuine dichotomy, since an unconditional basis
allows $W$ to be decomposed in uncountably many different ways. Note
also that the last two results imply that a homogeneous space $X$ must
have an unconditional basis, since otherwise no subspace of $X$ has
an unconditional basis, which implies, by Theorem 1.4, that some
subspace of $X$ is hereditarily indecomposable, which implies, by
homogeneity, that $X$ is hereditarily indecomposable, which implies,
by Theorem 1.3, that $X$ is not homogeneous. However, we have already
seen that a homogeneous space with an unconditional basis must be
isomorphic to $\ell_2$, so Banach's problem is solved once Theorem
1.4 is established.
\bigskip

The definitions we have given of an unconditional basis and a
hereditarily indecomposable Banach space are not the most convenient
for proving Theorem~1.4. Instead, we use simple characterizations
of these concepts in terms of finite sequences of vectors in the
space. First we need some more definitions and elementary facts 
concerned with bases. 

Because of Mazur's result that every Banach space has a subspace with
a basis, and indeed a basis with basis constant arbitrarily close to
1, we shall restrict attention to spaces with a monotone basis. Two
bases $(x_n)_{n=1}^\infty$ and $(y_n)_{n=1}^\infty$ are defined to be
$C$-{\it equivalent} if there exist constants $A$ and $B$ such that
$B/A\le C$ and
$$A\bnm{\sm n \infty a_nx_n}\le\bnm{\sm n \infty a_ny_n}
\le B\bnm{\sm n \infty a_nx_n}$$
for every sequence of scalars $(a_n)_{n=1}^\infty$ such that either
of the two sums converges. If $X$ is a Banach space with a given 
basis $(e_n)_{n=1}^\infty$, then the {\it support} of\break a vector $x=\sm
n \infty x_ne_n$, written $\supp(x)$, is the set of $n$ for which
$x_n\ne 0$.\break If $\max\supp(x)<\min\supp(y)$ we write $x<y$. A {\it block basis} of $X$ is defined to be a
sequence
$x_1<x_2<\dots$ of nonzero vectors. A subspace of $X$ generated by a block basis is
called a {\it block subspace}. The following very useful lemma of
Bessaga and Pelczynski (see [LT]) allows us to restrict our attention 
to block bases and block subspaces.

\proclaim{Lemma} Let $X$ be a Banach space with a basis
$(e_n)_{n=1}^\infty${\rm ,} let $Y$ be a subspace of $X$ and let
$\e>0$. Then $Y$ has a subspace $Z$ generated by a basis which is
$(1+\e)$\/{\rm -}\/equivalent to a block basis of $(e_n)_{n=1}^\infty$.
\endproclaim

In fact, more is true. For every sequence $(\d_n)_{n=1}^\infty$
of positive real numbers, we can find a subspace $Z$ with a
basis $(z_n)_{n=1}^\infty$ such that there is a normalized block basis
$(x_n)_{n=1}^\infty$ of $(e_n)_{n=1}^\infty$ with $\nm{x_n-z_n}\le\d_n$
for every $n$. In other words, $Z$ can be chosen to be an arbitrarily
small perturbation of a block subspace of~$X$. Moreover, given any
basic sequence in $Y$, we can choose $Z$ to be spanned by a subsequence.

We now state two simple lemmas. Proofs can be found in [LT] and [GM1].
From now on, when we use block-basis terminology, we shall assume
implicitly that the Banach space under discussion comes with some
particular chosen normalized monotone basis $(e_n)_{n=1}^\infty$. This
will save a good deal of writing.

\proclaim{Lemma}  Let $X$ be a Banach space. The following two
statements are equivalent\/{\rm :}
\begin{itemize}
\ritem{(i)} $X$ has no subspace with an unconditional basis.
\ritem{(ii)} For every block subspace $Y$ of $X$ and every
real number $C$ there is a sequence $\sleq y n$ of vectors in $Y$
such that
$$\bnm{\sm i ny_i}>C\bnm{\sm i n(-1)^iy_i}\ .$$
\end{itemize}

\endproclaim

\proclaim{Lemma}  Let $X$ be a Banach space. The following 
two statements are equivalent\/{\rm :}
\begin{itemize}
\ritem{(i)} $X$ is hereditarily indecomposable.
\ritem{(ii)} For every pair of block subspaces $Y,Z$ of $X$
and every real number $C$ there is a sequence $y_1<z_1<y_2<z_2<
\dots<y_n<z_n$ of vectors such that $y_i\in Y$ and $z_i\in Z$
for every $i${\rm ,} and such that
$$\bnm{\sm i n(y_i+z_i)}>C\bnm{\sm i n(y_i-z_i)}\ .$$
\end{itemize}

\endproclaim

Notice that if we insist that $Y=Z$ in condition (ii) of Lemma 1.7 then
we recover condition (ii) of Lemma 1.6 (with $n$ replaced by
$2n$). This makes it clearer what we need to do to prove Theorem
1.4. Let us say that a sequence $\sleq y n$ is $C$-{\it conditional}
if
$$\bnm{\sm i ny_i}>C\bnm{\sm i n(-1)^iy_i}\ .$$
Then we begin with a space $X$ such that every block subspace
contains a\break $C$-conditional sequence for every $C$, and must
find a (block) subspace with the stronger property that for
any {\it two} block subspaces and any $C$, a $C$-conditional
sequence can be found with its odd terms in one of the subspaces
and its even terms in the other. As we shall see later in the
paper, a far more general result is true, one which does not
use in a strong way the definition of $C$-conditional.

\section{Further definitions and some preliminary results}
 
\def \Si{\Sigma}

In this section we shall state the first of our main Ramsey-theoretic
results, which implies Theorem 1.4 very easily. It will be proved in
the next section. Given a Banach space $X$ (with a specified basis)
let us define $\Sigma_f=\Si_f(X)$ to be the set of all finite sequences
$\sleq x n$ of nonzero vectors in the unit ball of $X$.  Given an
arbitrary subset $\s\subset\Sigma_f$ (in other words, given some set of
finite block bases) the following two-player game can be defined,
between players S and P.  On the $n^{\rm th}$ move of the game, player
S chooses a block subspace $X_n\subset X$ (infinite-dimensional), and
player P chooses some point $x_n\in X_n$. The aim of P is to construct
a sequence $(\seq x n)\in\s$, whereas the aim of S is that at no stage
should the sequence constructed by P belong to $\s$.

Formally, a {\it strategy for} P is a function $\phi$ which, for any
finite block basis $\seq x n$ and any subspace $Y\subset X$, gives a
vector $x=\phi(\seq x n;Y)\in Y$. The strategy $\phi$ is a {\it winning strategy for} P if, given any sequence $X_1,X_2,\dots$ of
subspaces of $X$, the sequence $(x_1,x_2,\dots)$ defined inductively
by $x_1=\phi(\emptyset;X_1)$ and $x_{n+1}=\phi(\seq x n;X_{n+1})$ is
in $\s$.

If $\D=(\d_1,\d_2,\dots)$ is a sequence of positive scalars and
$\s\subset\Sigma_f$, then the $\D$-{\it expansion} of $\s$, denoted
$\s_\D$, is the set of block bases $(\seq x n)\in\Sigma_f$ such that
there exists $(\seq y n)\in\s$ with $\nm{y_i-x_i}\le\d_i$ for
every $i$. We shall also define $\s_{-\D}$ to be $(((\s)^c)_\D)^c$, or
in other words the set of block bases $(\seq x n)\in\Sigma_f$ such
that every $(\seq y n)\in\s$ with $\nm{y_i-x_i}\le\d_i$ for all
$i$ is in $\s$. It is easy to check that 
$(\s_{-\D})_\D\subset\s\subset(\s_\D)_{-\D}$.

\def \sq #1{(#1_n)_{n=1}^\infty}
\def \sqf #1 #2{(#1_n)_{n=1}^#2}

We now introduce some notation. If $A=(\seq y m)$ and $Y\subset X$ is
a subspace, then $[A;Y]$ will stand for the set of sequences $\sqf z
N\in\Sigma_f$ such that $z_i=y_i$ for $i\le m$ and $z_i\in Y$ for
$m<i\le N$. If $A$ is the null sequence then we shall write
$[Y]$. Given $\s\subset\Sigma_f$, we write $\s[A;Y]$ for the set of
sequences $\sqf x N$ such that $(\seq y m,\seq x N)\in
[A;Y]\cap\s$. Finally, when we say that P {\it has a winning strategy
for the game} $\s[A;Y]$, S's moves are understood to have to be
subspaces of $Y$. 

Our first Ramsey-type theorem is the following.

\proclaim{Theorem}   Let $X$ be a Banach space{\rm ,} let $\sigma$
be any subset of $\Si_f(X)$ and let $\D$ be a sequence of positive
real numbers. Then $X$ has a subspace $Y$ such that either
$\s[Y]=\emptyset$ or $P$ has a winning strategy for the game 
$\s_\D[Y]$. 
\endproclaim

 We refer to this as a Ramsey theorem for two reasons.
The first is that its proof closely resembles existing arguments in
Ramsey theory. However, even the statement can be regarded as
Ramsey-theoretical. If we call sequences in $\s$ blue and those
not in $\s$ red, then the theorem gives us a subspace such that
either every finite block sequence is red or there is such an
abundance of small perturbations of blue sequences that P has
a winning strategy for obtaining them. One might hope to find
a subspace where {\it all} sequences were close to blue sequences,
but a strengthened statement along these lines is false (for
nontrivial reasons, see the appendix).

Many arguments in infinite Ramsey theory depend on one particular
diagonalization procedure. It will be used often enough in this paper
for it to be well worth stating as an abstract principle.  We shall
need some more notation. A $*$-{\it pair} is a pair $(A,Z)$ where
$A\in\Si_f$, $Z$ is an infinite-dimensional block subspace and $x<z$ for
every $x\in A$ and $z\in Z$. We shall sometimes refer to $*$-pairs
simply as pairs, and we shall use the notation $A<Z$ as shorthand for
the support condition above. If $\D$ is a sequence $(\d_1,\d_2,\dots)$
of positive reals and $A=(\seq x n)$ and $B=(\seq y n)$ are block
bases of the same size, we use the shorthand notation $d(A,B)\le\D$ to
mean $d(x_i,y_i)\le\d_i$ for each $i$. The definition of a $\D$-{\it net} of a set of block bases is obvious. If $\Pi$ is a set of
$*$-pairs, we write $\Pi_\D$ for the set of pairs $(A,Z)$ such that
there is a block basis $B$ with $d(A,B)\le\D$ and $(B,Z)\in\Pi$.  If
$\sleq x n$, then we shall write $\sp{\seq x n}$ for the subspace
generated by $\seq x n$. Given a sequence $(\seq x n)\in\Si$, we
shall write $\Si_f(\seq x n)$ for the set of all sequences 
$(\seq y k)\in\Si_f$ such that every $y_i$ belongs to the subspace
$\sp{\seq x n}$.
 
\def \dim{\mathop{\rm dim}}

\proclaim{Lemma} Let $\D_1,\D_2,\dots$ be a sequence of
positive real sequences and let $\Pi_1,\Pi_2,\dots$ be a sequence of
sets of $*$-pairs satisfying the following conditions\/{\rm :}\/
\begin{itemize} 
\ritem{(i)} for every pair $(A,Z)$ and every $n$ there is $Z'\subset Z$ such that
$(A,Z')\in\Pi_n${\rm ;}

\ritem{(ii)} if $(A,Z)\in\Pi_n$ and $Z'\subset Z$ then $(A,Z')\in\Pi_n$.
\end{itemize}
 Then there exists a subspace $Y\subset X$ such that
$(A,Z)\in(\Pi_n)_{\D_n}$ for every pair $(A,Z)$ such that
$A$ is of length at least $n$ and both $A$ and $Z$ are subsets of $Y$. 
\endproclaim

\demo{Proof} Choose a block basis $y_1,y_2,\dots$ and a sequence of block
subspaces $X=Y_0\supset Y_1\supset Y_2\supset\dots$ inductively as
follows. Once we have chosen $\seq y {{n-1}}$ and $\seq Y {{n-1}}$,
let $y_n\in Y_{n-1}$ be arbitrary (except that\break $(\seq y n)\in\Si_f$) 
and let $\seq A N$ be a $\D_n$-net of the set $\Si_f(\seq y n)$.

By property (i) we can pick a sequence of subspaces $$Y_{n-1}\supset
V_{11}\supset\dots\supset V_{1n}\supset V_{21}\supset\dots\supset
V_{2n}\supset\dots\supset V_{N1}\supset\dots\supset V_{Nn}$$ such that
$(A_i,V_{ij})\in\Pi_j$ for every $1\le i\le N$ and $1\le j\le n$. Let
$Y_n=V_{Nn}$. By property (ii) we have $(A_i,Z)\in\Pi_j$ for every
$i,j$ and every $Z\subset Y_n$. Since $\sqf A N$ is a $\D_n$-net, we
find that $(A,Z)\in(\Pi_j)_{\D_n}$ for every $j\le n$ and every pair
$(A,Z)$ such that $A\in\Si(\seq y n)$ and $Z\subset Y_n$.

Let $Y=\sp{y_1,y_2,\dots}$ and suppose that $(A,Z)$ is a $*$-pair
inside $Y$, with $A$ of length $m\ge n$. Then there exists $k\ge n$
such that $A\subset\sp{\seq y k}$ and
$Z\subset\sp{y_{k+1},y_{k+2},\dots}\subset Y_k$. By our construction
therefore, $(A,Z)\in(\Pi_n)_{\D_n}$. Hence, $Y$ will do.
\enddemo

We shall often use specializations of Lemma 2.2. It may be convenient
for the reader if we state them separately. Let us define a {\it
singleton} $*$-{\it pair} to be a pair $(x,Z)$ where $x$ is a nonzero
vector of norm at most 1, $Z$ is a block subspace and $x<z$ for every
$z\in Z$ (or in other words a $*$-pair $(A,Z)$ for which $A$ is a
singleton). Given $\d>0$ and a set $\Pi$ of singleton $*$-pairs, write
$\Pi_\d$ for the set of $*$-pairs $(x,Z)$ for which there exists $x'$
such that $d(x,x')\le\d$, $x'<z$ for every $z\in Z$ and
$(x',Z)\in\Pi$.

\proclaim{{C}orollary} Let $\d>0$ and let $\Pi$ be a 
set of singleton $*$-pairs satisfying the following conditions\/{\rm :}
\begin{itemize}    
\ritem{(i)} for every pair $(y,Z)$ there is $Z'\subset Z$ such that
$(y,Z')\in\Pi$\/{\rm ;}\/

\ritem{(ii)} if $(y,Z)\in\Pi$ and $Z'\subset Z$ then $(y,Z')\in\Pi$.
\end{itemize}
 Then there exists a subspace $Y\subset X$ such that
every pair $(y,Z)$ with $y\in Y$ and $Z\subset Y$ belongs
to $\Pi_\d$. 
\endproclaim

\demo{Proof} Apply Lemma 2.2 with the following choices for the $\Pi_i$
and $\D_i$.  Let $\Pi_1$ be the set of all $*$-pairs $(A,Z)$ such that
if $A$ is a singleton, then $(A,Z)\in\Pi$. (Thus, every pair $(A,Z)$
for which $A$ is not a singleton belongs to $\Pi_1$.)  If $i>1$, then
let $\Pi_i$ be the set of all $*$-pairs. Let all $\D_i$ be the sequence
$(\d,1,1,1,\dots)$ (but all that matters is that the first term of
$\D_1$ should be $\d$). \enddemo

\proclaim{{C}orollary} Let $\D_1,\D_2,\dots$ and
$\Pi_1,\Pi_2,\dots$ be as in Lemma {\rm 2.2} but satisfying the following
additional condition\/{\rm :}
\begin{itemize}
\ritem{(iii)} if $(A,Z)\in\Pi_n$ then $A$ has length $n$.
\end{itemize}
 Then there exists a subspace $Y\subset X$ such that
$(A,Z)\in(\Pi_n)_{\D_n}$ for every pair $(A,Z)$ such that $A$ and
$Z$ are subsets of $Y$ and $A$ has length $n$. 
\endproclaim

{\it Proof}.  Apply Lemma 2.2 replacing each $\Pi_n$ by $\Pi_n'$, where
$\Pi_n'$ is the set of $*$-pairs $(A,Z)$ such that if $A$ has
length $n$ then $(A,Z)\in\Pi_n$. (As in Corollary~2.3, if
$A$ has any other length then $(A,Z)$ belongs to $\Pi_n'$.)
\hfill\qed

\section{A Banach-space dichotomy} 

In this section we shall prove Theorem 2.1 and show in detail how it
implies Theorem 1.4. This will complete the solution of Banach's
problem on homogeneous spaces. Before stating the next result, let us
make two more definitions. If $X$ is a Banach space, $Y$ is a block
subspace of $X$ and $\s\subset\Sigma_f(X)$, we shall say that $\s$ is
{\it large for} $Y$ if every block subspace of $Y$ contains a sequence
in $\s$, and {\it strategically large for} $Y$ if P has a winning
strategy for the game $\s[Y]$. More generally, if $(A,Y)$ is a
$*$-pair, we say that $\s$ is {\it large for} $[A;Y]$ if every block
subspace of $Y$ contains a sequence $B\in\Sigma_f$ such that
$(A,B)\in\sigma$, and {\it strategically large for} $[A;Y]$ if P has a
winning strategy for the game $\s[A;Y]$. Note that $\s$ is
(strategically) large for $[A;Y]$ if and only if $\s[A;Y]$ is
(strategically) large for $Y$.

\proclaim{Theorem}  Let $X$ be a Banach space with a given
monotone basis. Let $\Theta=\sq\theta$ and $\Delta=\sq\d$ be sequences
of positive real numbers such that $2\sum_{i=N}^\infty\d_i\le\theta_N$
for every $N$. If $\sigma_{-\Theta}$ is large for $X${\rm ,} then $X$ has a
block subspace $Y$ such that $\s_{2\D}$ is strategically large for
$Y$.
\endproclaim
 
\demo{Proof} Suppose that $\s\subset\Sigma_f$ is a set for which the result 
is false. Then $\sigma_{-\Theta}$ is large for $X$, so in particular
$\s$ is large for $X$; on the other hand, $\s_{2\D}$ is not
strategically large for any subspace of $X$. Let $\rho$ be the set of
sequences $(\seq x n)\in\s$ such that if $\sleq y k$ and $\sp{\seq y
k}$ is a proper subspace of $\sp{\seq x n}$, then $(\seq y
k)\notin\s$. It is easy to see that $\rho$ is still large for $X$ and
that $\rho_{2\D}$ is not strategically large for any subspace of $X$.

For each $n\ge 0$ Let $\D_n=(\seq \d n,0,0,\dots)$ and let
$\G_n=2\D-\D_n=(\seq \d n,2\d_{n+1},2\d_{n+2},\dots)$. We
now construct sequences $x_1,x_2,\dots$ and $X=X_0\supset X_1\supset
X_2\supset\dots$ with the following properties for every $n$:
\vglue4pt
\vglue2pt{(i)} $x_n\in X_{n-1}$;
\vglue2pt{(ii)} $\rho_{\D_n}$ is large for $[\seq x n;X_n]$;
\vglue4pt{(iii)} $\rho_{\G_n}$ is not strategically large for any $[\seq
x n;Z]$ with $Z\subset X_n$.  
\vglue4pt

The induction starts with the space $X_0$, since $\rho_{\D_0}$ is 
large for $X_0$ but $\rho_{\G_0}$ is not large for any subspace
of $X_0$. Having found $\seq x n$ and $\seq X n$, suppose we cannot find
suitable candidates for $x_{n+1}$ and $X_{n+1}$. Then for every 
$x\in X_n$ and every subspace $Y$ of $X_n$ we can find a subspace
$Z\subset Y$ such that either
\vglue4pt
{(a)} $\rho_{\D_{n+1}}\cap[\seq x n,x;Z]$ is empty
\vglue4pt
\noindent or
\vglue4pt {(b)} $\rho_{\G_{n+1}}$ is strategically large for
$[\seq x n,x;Z]$.
\vglue4pt

Let $\Pi$ be the set of singleton $*$-pairs $(x,Z)$ with $x\in X_n$
and $Z\subset X_n$ such that either (a) or (b) holds. We shall apply
Corollary 2.3. Condition (ii) of this corollary is obvious, and we
have just shown that condition (i) holds as well. Applying the
corollary with $\d=\d_{n+1}$, we obtain a subspace $Y\subset X_{n+1}$
such that, for every $*$-pair $(y,Z)$ with $y\in Y$ and $Z\subset Y$,
there exists $x$ with $\nm{y-x}\le\d_{n+1}$ such that either (a) or
(b) holds for the pair $(x,Z)$. Since $\nm{y-x}\le\d_{n+1}$, the first
alternative implies that $\rho_{\D_n}\cap[\seq x n,y;Z]$ is empty and
the second implies that $\rho_{\G_n}$ is strategically large for $[\seq
x n,y;Z]$. In particular, one of these two conclusions is true when
$Z=Y$. (There is a small technical point which is that $(y,Y)$ is not
a $*$-pair. However, the conclusion can be seen by considering for
each $y\in Y$ the $*$-pair $(y,Z)$, where $Z=\{z\in Y:y<z\}$.)

The set of $y$ such that $\rho_{\D_n}\cap[\seq x n,y;Y]$ is empty
cannot contain a subspace $Z$ of $Y$ since then 
$\rho_{\D_n}\cap[\seq x n;Z]$ is empty, which contradicts (ii) of
our inductive hypothesis. Therefore, the set of $y$ such that
$\rho_{\G_n}$ is strategically large for $[\seq x n,y;Y]$ is 
large for $Y$. But this gives P a winning strategy for the
game $\rho_{\G_n}[\seq x n;Y]$ and contradicts (iii) of our 
inductive hypothesis.

We now claim that the subspace generated by $\sq x$ has empty
intersection with $\sigma_{-\Theta}$. Indeed, let $(\seq z k)\in\Sigma_f$
be a sequence contained in the subspace $\sp{\seq x n}$. By property
(ii), the sequence $(\seq x n)$ can be extended to a sequence $(\seq x
m)$ in $\rho_{\D}$. Choose $(x_1',\dots,x_m')\in\rho$ such that
$\nm{x_i-x_i'}\le \d_i$ for every $i\le m$. Let $(z_1',\dots,z_k')$
be the corresponding perturbation of $(\seq z k)$, and notice that the
minimality condition satisfied by $\rho$ ensures that
$(z_1',\dots,z_k')\notin\sigma$. Our choice of $\D$ guarantees that
$\nm{z_j-z_j'}\le\Theta_j$ for every $j\le k$.  But this proves that
$(\seq z k)\notin\sigma_{-\Theta}$ and our claim is proved. This
contradicts the assumption that $\sigma_{-\Theta}$ was large.
\enddemo

 It is not hard to see that the above result is equivalent
to Theorem 2.1. 

\demo{Proof of Theorem {\rm 2.1}}  Let $\tau=\sigma_{\D/2}$. 
The assumption of the theorem and the fact that
$\sigma\subset\tau_{-\D/2}$ tell us that $\tau_{-\D/2}$ is large for
$X$. Therefore, applying Theorem 3.1 to $\tau$ (with $\Theta$ replaced
by $\D/2$) we can find some positive real sequence $\G\le\D/2$ such
that $\tau_\G$ is strategically large for some subspace of $X$. But
$\tau_\G\subset\s_\D$, so the result is proved. \pagebreak
\enddemo

Before we apply Theorem 2.1, let us introduce two further definitions.
We shall say that a block basis $\sq x$ is $C$-{\it unconditional} if
it generates a subspace which contains no $C$-conditional finite
sequences of blocks. ($C$-conditional sequences were defined just
after Lemma 1.7.)  We shall say that a Banach space $X$ is $C$-{\it
hereditarily indecomposable} if for every pair $Y,Z$ of block
subspaces of $X$ we can find $y\in Y$ and $z\in Z$ such that
$\nm{y+z}>C\nm{y-z}$.  This is equivalent to condition (ii) of Lemma
1.7 for that particular value of $C$. Lemma 1.7 asserts that $X$ is
hereditarily indecomposable if and only if it is $C$-hereditarily
indecomposable for every $C$.

\proclaim{{C}orollary}  Let $X$ be a Banach space. Then either
$X$ contains a $C$\/{\rm -}\/unconditional block basis or for every $\e>0$
it has a $(C-\e)$\/{\rm -}\/hereditarily indecomposable block subspace.
\endproclaim

\demo{Proof} Let $\s$ be the set of all sequences $(\seq x n)\in\Sigma_f$ 
that are\break $C$-conditional and contain at least one vector $x_i$ of norm
1. If $X$ contains no $C$-unconditional block basis, then $\s$ is
large, since every subspace contains a $C$-conditional sequence $(\seq
y n)$, and if we divide this by the largest value of $\nm{y_i}$ we
obtain a sequence in $\s$. Therefore, by Theorem 2.1, for any $\D>0$
we can find a block subspace $W$ of $X$ such that $\s_\D$ is
strategically large for $W$. Let us choose $\D$ such that
$\sm i \infty\d_i=\eta$ for some $\eta>0$ satisfying the inequality 
$(1+2\eta)^{-1}(C-2\eta)\ge C-\e$.

Now let $Y$ and $Z$ be arbitrary block subspaces of $W$ and consider
the strategy for S, the subspace player, of alternating $Y$ and $Z$.
Since $\s_\D$ is strategically large for $W$, P can defeat this
strategy, which means that P can choose a sequence
$(y_1,z_1,y_2,z_2,\dots,y_n,z_n)\in\s_\D$ such that $y_i$ belongs to
$Y$ and $z_i$ belongs to $Z$ for every $i$. (An unimportant technical
point is that P's strategy may succeed after an odd number of moves.
One can either alter the statement of Lemma 1.7 or let P choose a
sufficiently small $z_n$ to finish. The second approach is possible
because of the strict inequality in the definition of $C$-conditional
sequences.) We can then find $(y_1',z_1',\dots,y_n',z_n')\in\s$ such
that $\nm{y_i-y_i'}\le\d_{2i-1}$ and $\nm{z_i-z_i'}\le\d_{2i}$ for
every $i\le n$. Since the basis of $X$ is monotone, the norms
$\bnm{\sm i n(y_i'+z_i')}$ and $\bnm{\sm i n(y_i'-z_i')}$ are both at
least 1/2.  We therefore know that
$$\bnm{\sm i n(y_i'+z_i')}>C\bnm{\sm i n(y_i'-z_i')}\ge C/2\ . $$
By the triangle inequality and our choice of $\D$, we know that
$$\bnm{\sm i n(y_i+z_i)}\ge\bnm{\sm i n(y_i'+z_i')}-\eta$$ 
and 
$$\bnm{\sm i n(y_i-z_i)}\le\bnm{\sm i n(y_i'-z_i')}+\eta\ .$$ 
It follows from our choice of $\eta$ that 
$$\bnm{\sm i n(y_i+z_i)}>(C-\e)\bnm{\sm i n(y_i-z_i)}\ .$$
Since $Y$ and $Z$ were arbitrary subspaces of $W$, we have shown
that $W$ is $(C-\e)$-hereditarily indecomposable.
\enddemo

 A simple diagonalization now completes the proof of Theorem 1.4. 
 
\demo{Proof of Theorem {\rm 1.4}}  If $X$ has no subspace
with an unconditional basis, then Lemma 1.6 tells us that, for every
$C$, every block subspace of $X$ contains a $C$-conditional block
sequence. By Corollary 3.2 we can therefore find a nested sequence
$W_1\supset W_2\supset\dots$ of block subspaces of $X$ such that,
for every $n$, the subspace $W_n$ is $n$-hereditarily indecomposable.
Let $\sq w$ be a block basis of $X$ such that $w_n\in W_n$ for every
$n$ and let $W$ be the subspace generated by $\sq w$. We claim
that $W$ is hereditarily indecomposable. To see this, let $Y$
and $Z$ be arbitrary block subspaces of $W$ and let $C$ be any
real number. Choose a positive integer $n\ge C$. Then $Y\cap W_n$
and $Z\cap W_n$ are infinite-dimensional. Since $W_n$ is
$n$-hereditarily indecomposable, we can find a sequence satisfying
condition (ii) of Lemma 1.7 for this particular $C$. But $C$
was arbitrary, so the condition holds in general. Lemma 1.7 
therefore implies that $W$ is hereditarily indecomposable, as
claimed. \enddemo

\section{Definitions and preliminary results for 
infinite sequences}

We begin this section with a brief discussion of the connections
between our results so far and known results of infinite Ramsey
theory. The statement of Theorem 2.1 is strongly reminiscent of a
result in infinite Ramsey theory due to Nash-Williams [N-W], which
says the following. Let $\N^\om$ be the set of all infinite subsets of
$\N$. If $A$ is an open subset of $\N^\om$ (in the product topology)
then either $A$ or its complement contains all infinite subsets of
some $X\in\N^\om$. A set $A$ with this property is called a {\it Ramsey} set.

It is easy to check that Nash-Williams's result is equivalent to the
following statement. Let $\N^{<\om}$ be the set of all finite subsets
of $\N$ and let $A\subset\N^{<\om}$. Then there is an infinite subset
$X$ of $\N$ such that either no finite subset of $X$ is in $A$ or
for every infinite subset $Y$ of $X$ there exists $n$ such that 
$Y\cap\{1,2,\dots,n\}\in A$. Notice that if the first alternative
does not hold, then $A$ is {\it large} in an obvious sense. So
Nash-Williams's theorem asserts that if $A$ is a large subset of
$\N^{<\om}$ then there is an infinite subset $X$ of $\N$ for which
$A$ has a much stronger largeness property. This formulation
makes the resemblance with Theorem 2.1 very clear.

Nash-Williams's theorem was extended to all Borel sets by Galvin and
Prikry [GP]. A combinatorial lemma of theirs inspired our proof of
Theorem~2.1. Silver [S] proved that all analytic sets are Ramsey, and
Mathias proved that in a model constructed by Solovay all sets are
Ramsey. From these proofs there emerged a natural strengthening of the
Ramsey property, which was shown by Ellentuck [E] to be equivalent to
the property of Baire in a certain topology. For proofs of the
Galvin-Prikry lemma and Ellentuck's theorem, see [Bo] and for further
results in this direction, see [M].

The main aim of the next two sections is to extend Theorem 2.1 in a
similar way. However, as we will explain later (see the appendix), the
obvious analogue of Ellentuck's characterization is false, so we must
be satisfied with a result more like Silver's. In order to state it,
we shall need some definitions concerning infinite sequences and sets
of infinite sequences.

The obvious way to extend our earlier results would be to redefine
$\Si(X)$ as the set of all {\it infinite} sequences $x_1<x_2<\dots$ of
nonzero vectors in $X$ such that $\nm{x_n}\le 1$ for every $n$.
Instead, for technical reasons, our definition will be slightly
different: let $\Si(X)$ be the set of all sequences of pairs
$(x_1,\l_1),(x_2,\l_2),\dots$ where $x_1<x_2<\dots$ are vectors of
norm 1, and $\l_1,\l_2,\dots$ are real numbers in the interval
$[0,1]$. Of course, we can usually identify the pair $(x,\l)$ with the
vector $\l x$. The main difference between our definition and the
obvious definition is that if $x\ne y$, then we distinguish between
the pairs $(x,0)$ and $(y,0)$. In order to save writing, we shall
usually use a single letter to denote one of these pairs, unless it is
important to be careful. We shall refer to elements of $\Si$ as a
block bases. Sometimes we shall discuss finite block bases. Let us now
redefine $\Si_f$ to be the set of finite sequences
$\bigl((x_1,\l_1),\dots, (x_n,\l_n)\bigr)$ with the $(x_i,\l_i)$ as
above and $\sleq x n$. The {\it support} of a pair $(x_n,\l_n)$ is
defined to be the support of the vector $x_n$. The {\it subspace
generated by {\rm (}\/the pairs\/{\rm )}} $(x_1,\l_1),(x_2,\l_2),\dots$ is defined to
be the (block) subspace generated by the vectors $x_1,x_2,\dots\, $.
%\pagegoal=49pc

Several definitions to do with finite sequences can be easily adapted
for infinite sequences. For example, if $X$ is a Banach space, $Y$
is a block subspace of $X$ and $A=(\seq y m)$ is a finite sequence
of blocks (that is, pairs of the above kind) then $[A;Y]$ is now
defined as the set of all infinite sequences $\sq z\in\Si$ such
that $z_i=y_i$ for $i\le m$ and $z_i\in Y$ for $i>m$. Again we
denote this by $[Y]$ if $A$ is the null sequence. If $\s$ is a
subset of $\Si$, then we write $\s[A;Y]$ for the set of all sequences
$\sq w\in\Si(Y)$ such that $(\seq y m,w_1,w_2,\dots)\in\sigma$.
We say that $\s$ is {\it large for} $[A;Y]$ if every block subspace
of $Y$ contains a sequence in $\s[A;Y]$. Once again, $\s$ is large
for $[A;Y]$ if and only if $\s[A;Y]$ is large for~$Y$.

Given a subset $\s$ of $\Si$, we can define an infinite game just
as we did for finite sequences. The only difference is that P's
aim is to produce an infinite sequence that lies in $\s$, so the
game always lasts for ever (although one player may have a 
guaranteed win after finite time regardless of future moves).
We say that $\s$ is {\it strategically large for} $[A;Y]$ if 
P has a winning strategy for the set $\s[A;Y]$ when all of
S's moves are required to be subspaces of~$Y$.
%\pagegoal=48pc

Finally, we consider two topologies on $\Sigma$. Most of the time it
will be convenient to take as basic open sets all sets of the form
$\{(x_n)_{n=1}^\infty: x_n=y_n\ \hbox{for}\ 1\le n\le N\}$. In other
words, we put the discrete topology on $X$ and then take the product
topology on $\Sigma(X)$. However, for our main result we take a
different topology. First, define a metric on $S(X)\times[0,1]$ by
$d\bigl((x,\l),(y,\mu)\bigr)=\nm{x-y}+|\l-\mu|$. ($S(X)$ is the unit
sphere of $X$. Roughly speaking, we have made a hole where zero used
to be.) From this we derive a topology on $X$ and hence a different
product topology on $\Si(X)$. A basic open neighbourhood of the
sequence $((x_n,\l_n))_{n=1}^\infty$ is now a set of the form
$$\{((y_n,\mu_n))_{n=1}^\infty:\nm{y_i-x_i}+|\mu_i-\l_i|<\e_i\ 
\hbox{for every}\ i\le N\}$$
for some positive integer $N$ and positive real sequence $(\seq \e
N)$. The advantage of the first topology is that it is less messy to
talk about open and closed sets, and the advantage of the second is
that it makes $\Sigma(X)$ a Polish space (see Lemma 4.2 below), which
is more convenient for talking about analytic sets. (It is in order to
make $\Si$ complete metrizable that we allow ``zero vectors'' in a
block basis.)  However, a perturbation is involved in our result, so
the distinction between the two topologies is not at all important for
applications. We shall refer to ${\rm D}$-open and N-open sets and so on (for
``discrete'' and ``norm'') when it is not otherwise clear which
topology is meant.

We have not yet defined perturbations of infinite sequences. If $\s$
is a subset of $\Si(X)$ and $\D>0$ is an infinite sequence of positive
real numbers, then let $\s_\D$ denote the set of all sequences $\sq x$
such that there exists a sequence $\sq y\in\s$ with
$d(x_n,y_n)\le\d_n$ for every $n$. (Note that the $x_n$ and $y_n$ are
elements of $S(X)\times[0,1]$ and $d$ is the metric defined above.) We
now have enough notation to state the main theorem of this and the
next section.

\proclaim{Theorem}  Let $X$ be a Banach space{\rm ,} let 
$\s\subset\Si(X)$ be ${\rm N}$\/{\rm -}\/analytic and large for $X$ and let $\D>0$. Then
there is a subspace $Y$ of $X$ such that $\s_\D$ is strategically
large for $Y$.
\endproclaim

 Let us give a definition which will be useful for the
rest of the paper.

\demo{Definition} A set $\s\subset\Si(X)$ is {\it weakly Ramsey}
if for every $\D>0$ there is a subspace $Y\subset X$ such that
either $\s\cap[Y]$ is empty or $\s_\D$ is strategically large for $Y$.
\enddemo

 Thus, Theorem 4.1 states that ${\rm N}$-analytic sets are
weakly Ramsey, and a set $\s\subset\Sigma(X)$ is weakly Ramsey if the
conclusion of Theorem 4.1 holds for $\s$.

For the rest of this section we shall prove that ${\rm D}$-open sets and N-closed
sets are weakly Ramsey, and then we shall prove a few lemmas which are
needed for Theorem 4.1. First, we check that Theorem 2.1 still holds 
now that we have redefined $\Si_f$. One way of doing this is simply
to check that the proof we gave is still valid under the new 
interpretation. However, it may reassure the reader to see that one
can deduce the result for the new $\Si_f$ from Theorem~2.1.

\proclaim{Lemma} Let $X$ be a Banach space{\rm ,} let $\sigma$
be any subset of $\Si_f(X)$ and let $\D$ be a sequence of positive
real numbers. Then $X$ has a subspace $Y$ such that either
$\s[Y]=\emptyset$ or $P$ has a winning strategy for the game $\s_\D[Y]$.
\endproclaim

\demo{Proof} Define $\rho$ to be the set of all sequences
$$\bigl((1+\l_1)x_1/2,\dots,(1+\l_n)x_n/2\bigr)$$
such that the sequence
$$\bigl((x_1,\l_1),\dots,(x_n,\l_n)\bigr)$$
is in $\s$. It follows from the triangle inequality that
$$\nm{(1+\l)x-(1+\mu)y}\ge\max\{|\l-\mu|,\nm{x-y}-|\l-\mu|\}\ .$$
Therefore, applying Theorem 2.1 to $\rho$ with $\D$ replaced by $\D/4$  
we obtain the desired conclusion. \enddemo

\proclaim{Theorem}  {4.3}. All ${\rm D}$\/{\rm -}\/open sets are weakly Ramsey.
\endproclaim

\demo{Proof} This is just a reformulation of Lemma 4.2. Let $\s$ be a
${\rm D}$-open set. Let $\s_f$ be the set of finite sequences $(\seq x
n)\in\Si_f$ such that all their extensions to sequences in $\Si$ are
elements of $\s$. Since $\s$ is a ${\rm D}$-open set, every sequence in $\s$
has an initial segment belonging to $\s_f$. If $\s$ is large for $X$,
it follows that $\s_f$ is large for $X$ as well. Therefore, by Lemma
4.2, for every $\D>0$ there is a subspace $Y$ of $X$ such that
$(\s_f)_\D$ is strategically large for $Y$. If P uses this strategy
and plays for ever, then the resulting infinite sequence belongs to
$\s_\D$. Thus, $\s_\D$ is also strategically large for $Y$. This shows
that $\s$ is weakly Ramsey, as claimed.  \enddemo

Although it is weaker than our main result and not needed for the
proof, we shall now show that N-closed sets are weakly Ramsey, since
the proof is quite short, and we shall apply this weaker result
directly in \S 7.  As a matter of fact, we prove slightly more, by
defining a finer topology, which we shall call the D-$*$-topology, 
very similar to the so-called $*$-topology,\pagebreak or Mathias topology, on 
the infinite subsets of $\N$. We shall then show that D-$*$-closed
sets are close to being weakly Ramsey and that N-closed sets are
genuinely weakly Ramsey. In Section 7, we shall apply this result to
obtain a second Banach space dichotomy.

The D-$*$-topology is the topology generated by the sets $[A;Y]$
defined at the beginning of this section. Thus, a basic open
neighbourhood of a sequence $\sq x$ is some set of sequences $[A;Y]$,
where $A$ is an initial segment $(\seq x m)$ and $Y$ is a block
subspace containing $x_i$ for every $i>m$. (Similarly, one can define
the N-$*$-topology to have basic open sets of the following form: the
set of all sequences $\sq x\in\Sigma$ such that $\nm{x_i-y_i}<\d$ for
$i\le n$ and $x_i\in Y$ for $i>n$, where $\d>0$, $\sleq y n$ and $Y$
is an infinite-dimensional block subspace.) A set $\s$ is {\it completely weakly Ramsey} if, whenever $\s$ is large in $[A;Y]$ and
$\D>0$, there is a subspace $Z$ of $Y$ such that P has a winning
strategy for the game $\s_{\D}[A;Z]$. Note that if every analytic set
is weakly Ramsey, then every analytic set is completely weakly
Ramsey. We shall show in the appendix that the natural analogue of
Ellentuck's theorem, that a set is completely weakly Ramsey if and
only if it is a Baire set in the N-$*$-topology, is false. In fact, we
show that the intersection of two completely weakly Ramsey sets need
not be completely weakly Ramsey.

Given $A$ and $Y$ as above, we shall follow Galvin and Prikry by
saying that $Y$ {\it accepts} $A$ ({\it into} $\s$) if
$[A;Y]\subset\s$, and that $Y$ {\it rejects} $A$ ({\it from} $\s$) if
no subspace $Z\subset Y$ accepts $A$. Saying that $Y$ rejects the null
sequence from the complement of $\s$ is equivalent to saying that
every block subspace $Z\subset Y$ contains a sequence in $\s$, or in
other words that $\s$ is large in $Y$.

We are about to apply Corollary 2.4 several times. However, strictly
speaking one needs a different result concerning modified blocks
$(x,\l)$. Such a result can easily be proved in an identical way,
so we shall simply apply Corollary 2.4, interpreting the vectors
there as pairs $(x,\l)\in S(X)\times[0,1]$.

\proclaim{Lemma}  Let $\s\subset\Sigma$ be a set of infinite
sequences and let $\D_0>\D_1>\D_2>\dots$ be a sequence of positive
real sequences. Then there is a subspace $Y\subset X$ such that{\rm ,} for
any $k$ and any sequence $(\seq y k)\in\Si_f(Y)${\rm ,} either $Y$ accepts
$\seq y k$ into $\s_{\D_{k-1}}$ or $Y$ rejects $\seq y k$ from
$\s_{\D_k}$.
\endproclaim

\demo{Proof} For each $k$ let $\G_k={1\over 2}(\D_{k-1}+\D_k)$ and let
$\Theta_k={1\over 2}(\D_{k-1}-\D_k)$. Let $\Pi_k$ be the set of pairs
$(A,Z)$ such that $A$ is a sequence of length $k$ and $Z$ either
accepts $A$ into $\s_{\G_k}$ or rejects $A$ from $\s_{\G_k}$. It is
easy to see that the conditions of Corollary 2.4 hold, so we can pass
to a subspace such that, for every $k$ and every $A$ of size $k$,
every pair $(A,Z)$ is in $(\Pi_k)_{\Theta_k}$. Given any pair $(A,Z)$
with $A$ of length $k$, choose $A'\in\Si_f$ such that
$d(A,A')\le\Theta_k$ and $(A',Z)\in\Pi_k$. If $Z$ accepts $A'$ into
$\s_{\G_k}$, then $Z$ accepts $A$ into $\s_{\D_{k-1}}$ and if $Z$
rejects $A'$ from $\s_{\G_k}$ then $Z$ rejects $A$ from
$\s_{\D_k}$. This proves the result.
\enddemo

\proclaim{Theorem}   Let $\s$ be a {\rm D-}\/$*$\/{\rm -}\/closed subset of $\Sigma${\rm ,}
let $\tau$ be the complement of $\s$ and let $\D>0$. Then there exists
a subspace $Y\subset X$ such that either $P$ has a winning strategy
for the game $\s[Y]$ or $[Y]\subset\tau_\D$. 
\endproclaim

\demo{Proof} Suppose there is no subspace $Y$ such that $[Y]\subset\tau_\D$.  
Then $X$ rejects $\emptyset$ from $\tau_\D$. For $k=0,1,2,\dots$ let
$\D_k=2^{-k}\D$.  By Lemma 4.4, we may assume that, for every sequence
$(\seq x n)$, $X$ either accepts $(\seq x n)$ into $\tau_{\D_{n-1}}$
or rejects it from $\tau_{\D_n}$.

Now let P play the following strategy: for the $n^{{\rm th}}$ move,
choose $x_n$ such that $X$ rejects $(\seq x n)$ from $\tau_{\D_n}$.
We shall show both that this is possible and that it is a winning
strategy.

Suppose then that P has played the strategy successfully for $n$
moves.  We know that $X$ rejects $(\seq x n)$ from $\tau_{\D_n}$. Now
let S play the subspace $X_{n+1}$. If P cannot continue the strategy,
then, for every $x\in X_{n+1}$, $X$ does not reject $(\seq x n,x)$
from $\tau_{\D_{n+1}}$. But because of our assumption about $X$, this
implies that $X$ accepts $(\seq x n,x)$ into $\tau_{\D_n}$ for every
$x\in X_{n+1}$. It follows that $X_{n+1}$ accepts $(\seq x n)$ into
$\tau_{\D_n}$. But this contradicts the fact that $X$ rejects
$(\seq x n)$ from $\tau_{\D_n}$.

We have shown that the strategy is possible. Now let $(x_1,x_2,\dots)$
be a sequence produced by B when using it. Then we certainly know
that, for every $n$, $X$ rejects $(\seq x n)$ from $\tau$. That means
that, given any subspace $Y\subset X$, $[\seq x
n;Y]\cap\s\ne\emptyset$. Since $\s$ is $*$-closed, it follows that
$(x_1,x_2,\dots)\in\s$.
\enddemo

 Because of the way we defined weakly Ramsey sets, Theorem
4.5 does not imply that D-$*$-closed sets are weakly Ramsey, although
this is true up to an arbitrarily small perturbation. Moreover, if
$\s$ is D-$*$-closed, it does not follow that $\s_\D$ is D-$*$-closed.
Although it is not important for applications, we shall now show that
N-closed sets are weakly Ramsey, using the following simple result.

\proclaim{Lemma}   Let $\s$ be {\rm N-}\/closed and let $\D>0$. Then
$\s_\D$ is also {\rm N-}\/closed.
\endproclaim

\demo{Proof} Since $\s_\D=(\s_{\D/2})_{\D/2}$ we can assume that every
$\d_n$ is at most 1/2. We shall use bold face letters to denote
sequences in $\Si$. Let ${\bf x}$ belong to the closure of
$\s_\D$. Then there exist sequences ${\bf y}^m\in\s_\D$ converging
pointwise to ${\bf x}$. For every $m$ we can find a sequence ${\bf
z}^m\in\s$ such that $d({\bf y}^m,{\bf z}^m)\le\D$. For fixed $n$, the
vectors $z_n^m$ (or strictly speaking elements of $S(X)\times[0,1]$)
have bounded support, as otherwise $d(z_{n+1}^m,x_{n+1}^m)$ would be
at least $2/3$ for sufficiently large $m$. Therefore we can find a
sequence $m_1,m_2,\dots$ such that ${\bf z}^{m_i}$ converges
pointwise, to ${\bf z}$, say. Since $\s$ is N-closed, ${\bf z}$ must
be in $\s$. Because ${\bf z}^m\ra{\bf z}$ and ${\bf y}^m\ra {\bf x}$
pointwise we must have $d(z_n,x_n)\le\d_n$ for every $n$, which proves
that ${\bf x}\in\s_\D$, as required.
\enddemo

\proclaim{{C}orollary}  Every {\rm N-}\/closed set is weakly 
Ramsey.
\endproclaim

\demo{Proof} Let $\s$ be N-closed and let $\D>0$. Then $\s_{\D/2}$
is N-closed, by Lemma 4.6, and in particular D-$*$-closed.
Applying Theorem 4.5 to $\s_{\D/2}$ with $\D/2$ replacing $\D$
we obtain the result.
\enddemo

The rest of this section is devoted to lemmas which will be used to
prove our main result, Theorem 4.1.  If $\s\subset\Sigma$ has the
property that it is not large for any $[A;Y]$, then we shall say that
it is {\it nowhere large}, and otherwise that it is {\it large
somewhere}. In the usual infinite Ramsey theorems, it is shown that
being nowhere large is the same as being $*$-nowhere dense. That is
false here (see the appendix). We do, however, have the following
useful lemma, which states that the nowhere large sets are, up to
perturbation, closed under countable unions.

\proclaim{Lemma} Let $\s=\bigcup_{n=1}^\infty\s_n${\rm ,} let $\D_n>0$
for every $n$ and suppose that $(\s_n)_{\D_n}$ is nowhere large for
every $n$. Then $\s$ is nowhere large.
\endproclaim

\demo{Proof} Let $(B,W)$ be any pair. We must show that $\s$ is not large
for $[B;W]$. For each $n$ let $\Pi_n$ be the set of pairs $(A,Z)$ in
$W$ such that $A$ is a sequence of length exactly $n$ and
$(\s_n)_{\D_n}\cap[(B,A);Z]$ is empty. (We mean by $(B,A)$ the
concatenation of $A$ and $B$.) Our hypothesis implies that the sets
$\Pi_n$ satisfy the conditions for Corollary 2.4. We can therefore
find a subspace $Y\subset W$ such that $(A,Z)\in(\Pi_n)_{\D_n}$ for
every pair $(A,Z)$ in $Y$ with $A$ of length $n$. It is easy to see
that this implies that $[(B,A);Z]\cap\s_n=\emptyset$ for every such
pair. This says that no sequence $C$ of length $n$ in $Y$ can be added
to the end of $B$ and then extended in $Y$ to a sequence in $\s_n$. It
follows that $[B;Y]\cap\s=\emptyset$, because if we could find
$D\in\Si(Y)$ and $n$ such that $(B,D)\in\s_n$ then the first $n$ terms
of $D$ would give us a sequence $C$ contradicting what we have just
established. Since $B$ and $W$ were arbitrary, we have shown that $\s$
is nowhere large.  \enddemo

  For technical reasons we need a slight extension of Lemma 4.8
(Corollary 4.10 below).

\proclaim{Lemma}   If $\s$ is large for $[B;Y]$ and $\D>0${\rm ,} then for
every $m\in\N$ there exists a sequence $A$ of length $m$ in $Y$ and
a subspace $Z\subset Y$ such that $\s_\D$ is large in $[(B,A);Z]$.
\endproclaim

{\it Proof}. Suppose $\s$ is a counterexample to the assertion. Let $\Pi_m$ be
the set of pairs $(A,Z)$ in $Y$ such that $A$ is a sequence of length
$m$ and $\s_\D\cap[(B,A);Z]=\emptyset$. For $n\ne m$ let $\Pi_n$ be
the set of all pairs. Our assumption implies that the $\Pi_n$ satisfy
the conditions of Corollary 2.4. Therefore we can find a subspace
$W\subset Y$ such that $(A,Z)\in\Pi_m$ \pagebreak for every pair $(A,Z)$ in $W$
with $A$ of length $m$, which implies that $\s\cap[(B,A);Z]=\emptyset$
for every such pair. But this says that no sequence from $W$ of length
$m$ can be put on to the end of $B$ and extended in $W$ to a sequence
in $\s$. Therefore $\s$ is not large for $[B;Y]$.
\hfill\qed

\proclaim{{C}orollary}  If $\s$ is large somewhere{\rm ,}
$\s=\bigcup_{n=1}^\infty\s_n${\rm ,} $\D_1,\D_2,\dots$ are positive real
sequences and $k_1,k_2,\dots$ is a sequence of integers{\rm ,} then there
exists $n$ such that $(\s_n)_{\D_n}$ is large for some $[A;Z]$ with
$A$ a sequence of size at least~$k_n$.
\endproclaim

{\it Proof}. By Lemma 4.8 we can find $n$ and a pair $(B,Y)$ such that
$\s_{\D_n/2}$ is large for $[B;Y]$. Then by Lemma 4.9 with $m=k_n$,
$\s=\s_{\D_n/2}$ and $\D=\D_n/2$, we can find a sequence $C\in\Si_f$
of length $k_n$ and a subspace $Z\subset Y$ such that $\s_{\D_n}$ is
large in $[B,C;Z]$. The length of the sequence $(B,C)$ is certainly
at least $k_n$ so we are done.\hfill\qed\vglue3pt

 The next lemma is very similar in spirit to Lemma 4.4.

\proclaim{Lemma} Let $X$ be a Banach space{\rm ,} let
$\s_1,\s_2,\dots$ be a sequence of subsets of $\Si(X)$ and let
$\D_1,\D_2,\dots$ be a sequence of positive real sequences. Then there
exists a subspace $Y\subset X$ such that{\rm ,} for every $n$ and every
sequence $A\in\Si_f(Y)$ of length at least $n${\rm ,} either $\s_n\cap[A;Y]$
is empty or $(\s_n)_{\D_n}$ is large for $[A;Y]$.
\endproclaim

\demo{Proof} For each $n$ let $\Pi_n$ be the set of all $*$-pairs 
$(A,Z)$ such that either $(\s_n)_{\D_n/2}\cap[A;Z]$ is empty
or $(\s_n)_{\D_n/2}$ is large for $[A;Z]$. It is very easy to
check that the sets $\Pi_n$ satisfy the conditions of Lemma 2.2. 
Therefore, there is a subspace $Y$ such that every $*$-pair
$(A,Z)$ in $Y$ with $A$ of length at least $n$ belongs to
$(\Pi_n)_{\D_n/2}$. This says that there is some
$A'$ with $d(A,A')\le\D_n/2$ with $(A',Z)\in\Pi_n$. If
$(\s_n)_{\D_n/2}\cap[A',Z]=\emptyset$ then $\s_n\cap[A;Z]=
\emptyset$. If $(\s_n)_{\D_n/2}$ is large for $[A';Z]$,
then let $W$ be an arbitrary subset of $Z$. We can
choose an infinite sequence $B\in\Si(W)$ such that $A<B$
and $(A',B)\in(\s_n)_{\D_n/2}$. Then $(A,B)\in(\s_n)_{\D_n}$.
This shows that $(\s_n)_{\D_n}$ is large for $[A;Z]$. The
proof is complete.\enddemo

 We end this section with a lemma of a more technical 
nature.

\proclaim{Lemma} Let $\s\subset\Sigma$ be large{\rm ,} let
$\s=\bigcup_{n=1}^\infty\s_n${\rm ,} let $\D>0$ be a real sequence and let
$\psi$ be a function from $\N$ to $\N$. For any subspace $Y\subset X$
let $\rho(Y)$ be the set of sequences $B\in\Si_f(Y)$ such that
there exists $n$ with $\max\{n,\psi(n)\}$ at most the size of $B$
and $(\s_n)_\D$ large for $[B;Y]$. Then there exists a subspace
$Y\subset X$ such that $\rho(Y)$ is strategically large for $Y$.
\endproclaim

\demo{Proof} By Lemma 4.11 with each $\s_n$ replaced by $(\s_n)_{\D/3}$ and
each $\D_n$ equal to $\D/3$, we can find a subspace $Z$ such that, for
every $n\in\N$ and every sequence $A\in\Si_f(Z)$ of length at least
$n$, either $(\s_n)_{\D/3}\cap[A;Z]$ is empty or $(\s_n)_{2\D/3}$ is
large in $[A;Z]$.
 
Now $\s$ is large in $Z$ and therefore in any subspace of
$Z$. Therefore, by Corollary 4.10, for any subspace $W\subset Z$
there exist $n$, a subspace $V\subset W$ and a sequence
$B\in\Si_f(W)$ of length at least $\max\{n,\psi(n)\}$ such that
$(\s_n)_{\D/3}$ is large in $[B;V]$. In particular, $(\s_n)_{\D/3}
\cap[B;Z]$ is nonempty, so by our construction of $Z$ we know
that $(\s_n)_{2\D/3}$ is large for $[B;Z]$.

Define $\tau$ to be the set of sequences $B\in\Si_f(Z)$ such that
there exists $n$ with $\max\{n,\psi(n)\}$ at most the length of $B$
and $(\s_n)_{2\D/3}$ large for $[B;Z]$. We have just demonstrated
that $\tau$ is large for $Z$. By Theorem 2.1 there is a subspace 
$Y\subset Z$ such that $\tau_{\D/3}$ is strategically large for $Y$.
The result follows, since $\tau_{\D/3}\cap\Si_f(Y)\subset\rho(Y)$. 
\enddemo

\section{Analytic sets are weakly Ramsey}

In preparation for the proof of Theorem 4.1, we recall the definition of
and a basic fact about analytic sets. For a proof, and for general
background in descriptive set theory, we recommend Chapter 7 of Jech's
excellent book [J]. (This is the last chapter but the relevant part is
very clear and more or less independent of the rest of the book.)
 
\def \cn{{\cal N}}

The {\it Baire space} $\cal N$ is the set $\N^\N$ with the product
topology. Using continued fractions, one finds that $\cal N$ is
homeomorphic to the set of irrational numbers. If $X$ is a topological
space, an {\it analytic} subset of $X$ is the continuous image of a
Borel subset. If $X$ is a Polish space, this can be shown to be
equivalent to saying that the subset is a continuous image of $\cn$.
This equivalent definition is usually the easiest to use. It is clear
that every Borel set is analytic. Souslin's theorem (which we do not
need here) states that two disjoint analytic sets can be contained in
two disjoint Borel sets, which implies that a set is Borel if and only
if it and its complement are analytic.

 Before giving the proof of Theorem 4.1, let us quickly
verify that the set we shall be considering is indeed a Polish space.

\proclaim{Lemma} Let $X$ be a Banach space. Then $\Si(X)$ is
a Polish space in the {\rm N-}\/topology.
\endproclaim

\demo{Proof}  We must find a complete metric which gives rise to the
N-topology on $\Si(X)$. Such a metric is given by
$$d\Bigl((x_n,\l_n)_{n=1}^\infty,(y_n,\mu_n)_{n=1}^\infty\Bigr)
=\max n^{-1}d\bigl((x_n,\l_n),(y_n,\mu_n)\bigr)$$
where on the right-hand side we have taken the metric defined
earlier on $S(X)\times[0,1]$. Let $S_m$ be the sequence
$\bigl((x_{mn},\l_{mn})\bigr)_{n=1}^\infty$ and suppose that
the sequence $(S_m)_{m=1}^\infty$ is Cauchy in the above metric.
The only small point to note is that the supports of the vectors
$x_{mn}$ are bounded for fixed $n$, as otherwise the vectors
$x_{m,n+1}$ would not form a Cauchy sequence in $X$. It is therefore
easy to see that the sequences $S_m$ converge pointwise and hence
in the metric, and that the limit is in $\Si(X)$. It is also simple
to check that this metric gives rise to the N-topology on $\Si(X)$.
\enddemo

\demo{Proof of Theorem {\rm 4.1}} Since $\Sigma$ is a Polish space
(in the N-topology) we can use the equivalent definition of analytic
sets.  Suppose therefore that $\s\subset\Sigma$ and that there is a
continuous map $f:\cn\ra\Sigma$ such that $f(\cn)=\s$. Given a finite
sequence $\seq n k$ of integers, let $\s_{\seq n k}$ denote the image
under $f$ of the set of sequences in $\cn$ starting $\seq n k$. Note
that, for any $\D>0$, $(\s_{\seq n k})_\D
=\bigcup_{n_{k+1}=1}^\infty(\s_{\seq n {{k+1}}})_\D$.
\vglue4pt
Let $0=\D_0,\D_1,\D_2,\dots$ be a (pointwise) strictly increasing
sequence of real sequences less than $\D$. For $n\ge 0$ let $\G_n$
be the sequence ${1\over 4}(\D_{n+1}-\D_n)$ and for $i=1,2,3$
let $\D_{n,i}=\D_n+i\G_n$.

Let $\phi$ be a bijection between the set of finite sequences of
integers and~$\N$. By Lemma 4.11 we can find a subspace $W$ of $X$
such that, for every $(\seq n k)$ and every sequence $A\in\Si_f(W)$ of
length at least $\phi(\seq n k)$, either\break $(\s_{\seq n k})_{\D_{k,1}}
\cap[A;W]$ is empty or $(\s_{\seq n k})_{\D_{k,2}}$ is large for
$[A;W]$.

For every $(\seq n k)$ and every pair $(A,Z)$ in $W$, define
$\rho_{\seq n k}[A;Z]$ to be the set of all sequences
$B\in\Si_f(Z)$ such that there exists $n_{k+1}$ with
$\phi(\seq n {{k+1}})$ at most the length of $B$ and
$(\s_{\seq n {{k+1}}})_{\D_{k,3}}$ large for $[A,B;Z]$.
Let $\Pi_{\seq n k}$ be the set of pairs $(A,Z)$ such that
if $A$ has length at least $\phi(\seq n k)$ then either
\vglue4pt
{(a)} $(\s_{\seq n k})_{\D_{k,1}}\cap[A;Z]$ is empty
\vglue4pt
\noindent or
\vglue4pt {(b)} $\rho_{\seq n k}[A;Z]$ is strategically large 
for $Z$.
\vglue4pt
We shall apply Lemma 2.2 to the sets $\Pi_{\seq n k}$. However,
unlike with our previous diagonalization arguments, it will not
be trivial to verify that the relevant conditions are satisfied.
Instead, we use Lemma 4.12. Condition (ii) of Lemma 2.2 is of
course not a problem. As for condition (i), let $(A,Z)$ be
any pair and suppose that the length of $A$ is at least $\phi(\seq n k)$.
If (a) is not true for the pair $(A,Z)$, then by our construction
of the subspace $W$  we know that $(\s_{\seq n k})_{\D_{k,2}}$ is
large for $[A;W]$ and hence for $[A;Z]$. Now 
$$(\s_{\seq n k})_{\D_{k,2}}=\bigcup_{n_{k+1}=1}^\infty
(\s_{\seq n {{k+1}}})_{\D_{k,2}}\ .$$
Lemma 4.12 applied to this union with $\D=\G_k$ implies that there
exists a subspace $Z'\subset Z$ such that $\rho_{\seq n k}[A;Z']$ is
strategically large for $Z'$, which tells us that (b) is true of the
pair $(A,Z')$.  This shows that condition (ii) of Lemma 2.2 is indeed
satisfied.

Using that lemma, we can find a subspace $Y\subset W$ such that every
pair $(A,Z)$ in $Y$ with $A$ of length at least $\phi(\seq n k)$
belongs to $(\Pi_{\seq n k})_{\G_k}$.  Let us define perturbations of
the sets $\rho_{\seq n k}[A;Z]$ by defining $\tau_{\seq n k}[A;Z]$ to
be the set of all sequences $B\in\Si_f(Z)$ such that there exists
$n_{k+1}$ with $\phi(\seq n {{k+1}})$ at most the length of $B$ and
$(\s_{\seq n {{k+1}}})_{\D_{k+1}}$ large for $[A,B;Z]$. The 
statement that $(A,Z)$ belongs to $(\Pi_{\seq n k})_{\G_k}$ implies
that either\break $(\s_{\seq n k})_{\D_k}\cap[A;Z]$ is empty or
$\tau_{\seq n k}[A;Z]$ is strategically large for $Z$. This is
true in particular when $Z=Y$.

We claim now that $\s_\D$ is strategically large for the subspace
$Y$, and to prove it we shall describe a suitable strategy for P.
The game starts with $\s=\s_{\D_0}$ large for $Y$. By our construction
of $Y$, the set $\tau$ is strategically large for $Y$. That means
that P can ensure that after some finite number of moves (s)he 
will reach a sequence $B_1\in\Si_f(Y)$ such that there exists $n_1$
with $\phi(n_1)$ at most the length of $B_1$ and $(\s_{n_1})_{\D_1}$
large for $[B_1;Y]$. In general, suppose that P has reached a
sequence $B_k\in\Si_f(Y)$ such that there exists $(\seq n k)$
with $\phi(\seq n k)$ at most the length of $B_k$ and 
$(\s_{\seq n k})_{\D_k}$ large for $[B_k;Y]$. The construction
of $Y$ guarantees that $\tau_{\seq n k}$ is strategically large
for $Y$, which means that P can play in such a way as to extend
$B_k$ to a sequence $B_{k+1}$ such that there exists $n_{k+1}$
with $\phi(\seq n {{k+1}})$ at most the length of $B_{k+1}$
and $(\s_{\seq n {{k+1}}})_{\D_{k+1}}$ large for $[B_{k+1};Y]$.

Let P play as above and let $B$ be the infinite sequence which has the
finite sequences $B_k$ as initial segments. That is, $B$ is the
sequence produced by P at the end of the game. It remains to prove
that $B\in\s_\D$. We do this by showing that
$d(B,f(n_1,n_2\dots))\le\D$. Write $B=\sq x$ and $f(n_1,n_2,\dots)=\sq
y$. If our assertion is false then there exists $t$ such that
$d(x_t,y_t)>\d_t$. By the continuity of $f$, for every sufficiently
large $k$ and for every sequence $(m_1,m_2,\dots)\in\cn$ with
$m_i=n_i$ for all $i\le k$, if $z=f(m_1,m_2,\dots)$, then
$d(z_t,y_t)<\d_t-d(x_t,y_t)$.  For any $r\ge\max\{t,\phi(\seq n k)\}$,
we have just proved that $(\s_{\seq n k})_\D\cap[\seq x r;Y]$ is
empty, which implies that $(\s_{\seq n k})_{\D_k}$ is not large for
$[\seq x r;Y]$. Choosing $k$ such that $\phi(\seq n k)>t$, we find
that there does not exist $r\ge\phi(\seq n k)$ such that
$(\s_{\seq n k})_{\D_k}$ is large for $[B_k;Y]$, which is a
contradiction. \enddemo

\section{Spaces where the result can be strengthened}

We show first how to strengthen Theorem 4.1 to a genuine Ramsey result
when the Banach space given at the beginning is $c_0$ and the block
bases are normalized. The big difference between $c_0$ and other
spaces is the following theorem [G5].

\proclaim{Theorem}  Let $A_1\cup\dots\cup A_r$ be a partition of 
the unit sphere of $c_0$ and let $\e>0$. Then there exists an 
infinite\/{\rm -}\/dimensional block subspace $Y$ of $c_0$ and some $j\le r$
such that every point in the unit sphere of $Y$ is within $\e$ of 
some point of $A_j$.
\endproclaim

Define a subset $A$ of the unit sphere of $c_0$ to be
{\it large} if it has a nonempty intersection with every block
subspace of $c_0$. (The word {\it asymptotic} is more standard
in Banach space theory, but less in keeping with the terminology 
of this paper and other papers in Ramsey theory.) An equivalent
formulation of Theorem 1 is that if $A$ is a large subset of
$c_0$ and $\e>0$ then there exists an infinite-dimensional
subspace $Y$ of $c_0$ such that every point of $Y$ lies within
$\e$ of some point of $A$.

We shall now prove something similar for analytic sets of normalized
block bases of $c_0$ (or equivalently, of block subspaces of $c_0$).
Rather than giving full details, which are very similar to the proof
of Theorem 4.1, we merely indicate the points at which that argument
must be changed. Let $\Si_1(X)$ denote the set of all normalized block
bases of $X$ and let us now interpret notation such as $[X]$ to refer
to normalized block bases. (Alternatively, one could consider sets
$\s\subset\Si(X)$ such that $((x_n,\l_n))_{n=1}^\infty\in\s$ implies
that $((x_n,\mu_n))_{n=1}^\infty$ for all sequences $(\mu_n)_{n=1}^\infty$
of numbers in the interval $[0,1]$.) Let us define a subset $\s$ of
$\Si_1(c_0)$ to be {\it Ramsey} if for every $\D>0$ there is a block
subspace $Y$ of $c_0$ such that either $\s\cap[Y]$ is empty or
$[Y]\subset\s_\D$.

\proclaim{Lemma} {\rm N-}\/closed subsets of $\Si_1(c_0)$ are Ramsey.
\endproclaim

\demo{Proof} Corollary 4.7 implies that N-closed subsets of $\Si(c_0)$ are
weakly Ramsey. Let $\s$ be N-closed and large for $c_0$, let $\D>0$
and let $X$ be a block subspace of $c_0$ such that $\s_{\D/2}$ is
strategically large for $X$. Let
$\D/2=\D_0<\D_0^+<\D_1^-<\D_1<\D_1^+<\D_2^-<\D_2<\dots$ be a strictly
increasing sequence of sequences bounded above by $\D$. For each $k$
let $\Pi_k$ be the set of pairs $(A,Z)$ such that $A$ has length $k$
and either
\vglue4pt
 {(a)} for every $z\in Z$, $\s_{\D_k^+}$ is not strategically
large for $[(A,z);X]$

\noindent or

 {(b)} for every $z\in Z$, $\s_{\D_{k+1}^-}$ is strategically 
large for $[(A,z);X]$.
\vglue4pt
\noindent Clearly, if $(A,Z)\in\Pi_k$ and $Z'$ is a subspace of $Z$
then $(A,Z')\in\Pi_k$. We now verify condition (i) of Lemma 2.2.
Let $(A,Z)$ be a pair such that (a) is not true for any subspace
of $Z$. This says that the set of $z\in Z$ for which $\s_{\D_k^+}$
is strategically large for $[(A,z);X]$ is large for $Z$. Hence, by
Theorem 6.1 there is a subspace $Z'\subset Z$ such that for every
$z\in Z'$ the set $\s_{\D_{k+1}^-}$ is  strategically large for 
$[(A,z);X]$. (Note that we did not use the full strength of the
fact that $\D_k^+<\D_{k+1}^-$. The strict inequality was important
only for the $k^{\rm th}$ terms of these sequences.)

We may therefore apply Corollary 2.4, which, for a suitable choice
of perturbations, gives us a subspace $Y$ of $X$ such that, for
every finite normalized block basis $A$ in $Y$, either
\vglue4pt {(c)} for every $z\in Y$, $\s_{\D_k}$ is not strategically
large for $[(A,z);X]$
\vglue4pt
\noindent or
\vglue4pt {(d)} for every $z\in Y$, $\s_{\D_{k+1}}$ is strategically
large for $[(A,z);X]$.
\vglue4pt 
We claim now that $[Y]\subset\s_\D$. The proof is by induction. We
start with the statement that $\s_{\D_0}$ is strategically large for 
$[X]$. It follows that there exists $z\in S(Y)$ such that $\s_{\D_0}$
is strategically large for $X$ (or else S would have a winning
strategy for the game $\s_{\D_0}[X]$ by starting with the subspace 
$Y$). By construction of $[Y]$, for every $z_1\in S(Y)$ the set
$\s_{\D_1}$ is strategically large for $[z_1;X]$. By a similar
argument, for any fixed $z_1$ we then have that for every $z_2\in S(Y)$
the set $\s_{\D_2}$ is strategically large for $[z_1,z_2;X]$.
Continuing in this way we find that for any finite normalized
block basis $(\seq z k)$ in $Y$, $\s_{\D_k}$ 
is strategically large for $[\seq z k;X]$. It follows that
$(\seq z k)$ can be extended to a sequence in $\s_{\D}$. By
Lemma 4.6, $\s_{\D}$ is N-closed, and therefore D-closed.
It follows that every sequence $\sq z\in\Si_1(Y)$ belongs
to $\s_\D$.\enddemo

\proclaim{{C}orollary}  {\rm N-}\/open subsets of $\Si_1(c_0)$ are
Ramsey.
\endproclaim

{\it Proof}. Let $\s$ be an N-open subset of $\Si_1(c_0)$ and let
$\tau$ be the complement of $\s_\D$. It is not hard to show 
that $\s_\D$ is N-open, so $\tau$ is N-closed. If no subpace
$Y$ can be found such that $[Y]\subset\s_\D$, then $\tau$ is
large for $c_0$. Then, by Lemma 6.2, there is a subspace $Y$
such that $[Y]\subset\tau_\D$. However, $\tau_\D\cap\s$ is
empty, so $\s\cap[Y]$ is empty. \hfill\qed
\vglue4pt

  Given a subspace $Y$ of $c_0$ and a collection $\rho$
of finite sequences, let us say that $\rho$ is {\it very large for}
$Y$ if every sequence in $\Si_1(Y)$ has an initial segment in~$\rho$.
\vglue-6pt

\proclaim{Lemma} Let $\s\subset\Sigma_1(c_0)$ be large{\rm ,} let
$\s=\bigcup_{n=1}^\infty\s_n${\rm ,} let $\D>0$ be a real sequence and let
$\psi$ be a function from $\N$ to $\N$. For any subspace $Y\subset X$
let $\rho(Y)$ be the set of finite normalized block bases $B$ in $Y$
such that there exists $n$ with $\max\{n,\psi(n)\}$ at most the 
length of $B$ and $(\s_n)_\D$ large for $[B;Y]$. Then there exists
a subspace $Y\subset X$ such that $\rho(Y)$ is very large for $Y$. 
\endproclaim

{\it Proof}. For each $Y$, let $\tau(Y)$ be defined in the same way as
$\rho(Y)$ but with $\D$ replaced by $\D/2$. Lemma 4.12 implies that 
there is a subspace $W$ such that $\tau(W)$ is strategically large
for $W$. In particular, $\tau(W)$ is large for $W$. Hence, by
Corollary 6.3, there is a subspace $Y\subset W$ such that 
$[Y]\subset\tau(W)_{\D/2}$. The result now follows from the
fact that $\tau(W)_{\D/2}\subset\rho(W)$ and that $\rho(Y)\subset
[Y]\cap\rho(W)$. \phantom{Immad}
\hfill\qed 

\proclaim{Theorem}  Let $\s$ be an ${\rm N}$\/{\rm -}\/analytic subset of
$\Sigma_1(c_0)$ and let $\D>0$. Then there is a subspace $X\subset
c_0$ such that either $[X]\cap\s =\emptyset$ or $[X]\subset\s_\D$.
\endproclaim

\demo{{P}roof {\rm (}\/Sketch\/{\rm )}} The proof is very similar to that of
Theorem 4.1. The main changes are that the phrase ``strategically 
large'' should be replaced by ``very large'' and that instead of
Lemma 4.12 one should use Lemma 6.4. In the penultimate paragraph
of the proof, instead of considering a game between P and S, one
considers an arbitrary sequence $B$ in $\Si_1(Y)$. The construction
of $Y$ now guarantees that this sequence has an increasing 
collection of initial segments $B_1,B_2,\dots$ such that for
some sequence $n_1,n_2,\dots$ of positive integers, we have
for every $k$ that $\phi(\seq n k)$ is at most the length of
$B_k$ and $(\s_{\seq n k})_{\D_k}$ is large for $[B_k;Y]$.
The proof that the distance between $B$ and $f(n_1,n_2,\dots)$
is at most $\D$ is the same as in the proof of Theorem 4.1.~\hfill\qed
\enddemo

The above result also applies, for trivial reasons, to spaces
containing $c_0$. It follows from results of Milman [Mi] and Odell and
Schlumprecht [OS] that these are the only spaces for which this
stronger Ramsey result is true. The result can in fact be strengthened
further; it is basically symmetrical and so an Ellentuck-type theorem
holds. 

Going to the opposite extreme, let us consider the ``worst'' spaces and 
return to nonnormalized sequences. Recall that a space $X$ (which we
assume for convenience has a basis here) is {\it hereditarily
indecomposable} if it has one of the following equivalent properties:
\begin{itemize}
\ritem{(i)} $X$ has no subspace that can be decomposed as a
topological direct sum $Y\oplus Z$ with $Y$ and $Z$ 
infinite-dimensional;
\ritem{(ii)} for any pair of infinite-dimensional subspaces $Y$ and
$Z$ of $X$ and any $\e>0$ there exist $y\in Y$ and $z\in Z$ such
that $\e\nm{y+z}>\nm{y-z}$;
\ritem{(iii)} for any pair of infinite-dimensional block subspaces 
$Y$ and $Z$ of $X$ and any positive real sequence $\D$ there are 
block subspaces $Y'\subset Y$ and $Z'\subset Z$ such that 
$d(Y',Z')\le\D$.
\end{itemize}

The third property above, which we shall use, states that up to small
perturbations the subspaces of $X$ form a filter. The existence of
such spaces is proved in [GM1]. The observation we make here is that if
$X$ is hereditarily indecomposable and $\s\subset\Sigma(X)$ is a set
such that, for some subspace $Y\subset X$, P has a winning strategy
for the game $\s[Y]$, then for any $\D>0$ P has a winning strategy for
the game $\s_\D[X]$. Indeed, let $X_1,X_2,\dots$ be any sequence of
(infinite-dimensional block) subspaces of $X$. Using property (iii)
with the perturbations sufficiently small, we can find a sequence of
subspaces $Y_1,Y_2,\dots$ such that for every $n\in\N$ and every $y\in
Y_n$ of norm 1 there exists  \pagebreak $x\in X_n$ of norm 1 such that
$\nm{y-x}<\d_n$. If P plays a (fixed) winning strategy against the
moves $Y_1,Y_2,\dots$ from S, then P has produced a sequence in
$\s_\D$. Thus we have the following result.

\proclaim{Theorem}  Let $X$ be a hereditarily indecomposable space{\rm ,}
let $\D>0$ and let $\s\subset\Sigma(X)$ be ${\rm N}$\/{\rm -}\/analytic. Then either there
is a subspace $Y$ such that $\s\cap[Y]=\emptyset$ or there is a
winning strategy for $P$ for the game $\s_\D[X]$.
\endproclaim

\section{Applications}

Our aim in this section is to obtain the beginnings of a
classification of separable infinite-dimensional Banach spaces.  We do
not mean by this a complete description of all such spaces up to
isomorphism, since such a task is clearly hopeless. Rather, we wish to
find classes of Banach spaces with the following properties:
\begin{itemize}
\item[(a)] If $X$ is in one of the classes, then so is every subspace of $X$.

\item[(b)] Every space has a subspace in one of the classes.

\item[(c)] The classes are very obviously disjoint.

\item[(d)] Knowing that a space belongs to a particular class gives a lot of
information about the structure of the space and what operators can be
defined on it.
\end{itemize}

Properties (c) and (d) above are of course related. We begin with the
comment that it is slightly easier to deduce Theorem 1.4 from Theorem
4.1 than from Theorem 2.1, because the extra diagonalization (see the
end of \S 3) can be avoided. One defines $\s$ to be the set of
all infinite sequences $(x_n,\l_n)_{n=1}^\infty$ such that for every
positive integer $m$ there exist $r$ and $s$ such that the sequence
$(\l_rx_r,\l_{r+1}x_{r+1},\dots,\l_sx_s)$ is $m$-conditional. Lemma
1.6 implies that $X$ contains no unconditional basic sequence if and
only if $\s$ is large. It is easy to check that $\s$ is a ${\rm
G}_\d$-set. So if $X$ contains no unconditional basic sequence and
$\D>0$, then we can find a subspace $Y$ of $X$ such that $\s_\D$ is
strategically large for $Y$. For sufficiently small $\D$, this shows,
as in the proof of Corollary 3.2, that $Y$ satisfies condition (ii) of
Lemma 1.7, and is therefore hereditarily indecomposable.

In [GM1] it was shown that every operator on a (complex) hereditarily
indecomposable space is a strictly singular perturbation of a multiple
of the identity. This result was strengthened by Ferenczi [F2], who
proved that every operator from a subspace of a hereditarily
indecomposable space into the space itself is a strictly singular
perturbation of the inclusion map. We have therefore found two classes
satisfying properties (a) to (d). However, property (d) is a matter of
degree. It is not clear whether one would want to classify the
hereditarily indecomposable spaces any further (at least from the
point of view of understanding the operators defined on them) but
there is a wide variety of spaces with an unconditional basis
differing in important ways.  At one extreme there are the $\ell_p$
spaces, the ones with most symmetry, and at the other is a space
constructed in [G4] which is not isomorphic to any proper subspace of
itself. In fact, it is shown in [GM2] that every operator on the space
is a strictly singular perturbation of a diagonal map. (It is not
known whether every operator from a subspace to the whole space is a
strictly singular perturbation of the restriction of a diagonal map.)
Since the diagonal maps with bounded entries down the diagonal are
trivially continuous, such a space can be regarded as the ``worst''
type of space with an unconditional basis, just as a hereditarily
indecomposable space is the ``worst'' type of general space.

It is natural to ask, therefore, what can be said about a space with
an unconditional basis if it does not contain the ``worst'' kind of
subspace.  Is there some result stating that such a space must have a
``nice'' subspace? We shall give a result of this type, but we need
some definitions first.

Recall first that an infinite-dimensional Banach space $X$ is said to
be {\it minimal} if every infinite-dimensional subspace of $X$ has a
further subspace isomorphic to $X$. This term was coined by Rosenthal,
who conjectured that every infinite-dimensional space has a minimal
subspace. However, it was shown by Casazza and Odell [CO] that
Tsirelson's space gives a counterexample. Of course, a hereditarily
indecomposable space gives a much stronger counterexample, but for
this very reason Tsirelson's space is a more interesting one. It has
in particular a property which we shall now define and discuss in some
detail. Recall that two spaces $X$ and $Y$ are said to be {\it totally
incomparable} if no infinite-dimensional subspace of $X$ is isomorphic
to a subspace of $Y$. Let us say that a space $X$ is {\it quasi\/{\rm -}\/minimal} 
if it does not contain a pair of totally incomparable
subspaces. The reason for this terminology is given by the next simple
lemma. If $Y$ and $Z$ are two spaces, we write $Y\le Z$ to mean $Y$
embeds into $Z$.

\def \cs{{\cal S}}
 
\proclaim{Lemma} A space $X$ is quasi\/{\rm -}\/minimal if and
only if there is a collection $\cal S$ of subspaces of $X$ with the
following properties\/{\rm :}
\begin{itemize}
\ritem{(i)}  The set $\cal S$ is partially ordered by $\le$.

\ritem{(ii)}  If $Y,Z\in\cs$ then there exists $W\in\cs$ such
that $Y\ge W$ and $Z\ge W$.

\ritem{(iii)} Every subspace of $X$ has a subspace isomorphic to a 
subspace in $\cs$.  
\end{itemize} \pagebreak
\endproclaim

\demo{Proof} Write $Y\sim Z$ if $Y\le Z$ and $Z\le Y$. Then $\sim$ is
an equivalence relation. If $\cs$ consists of a representative from
each equivalence class then it clearly has the properties stated.
The converse is obvious from properties (ii) and (iii).
\enddemo

{\it Remarks}. The important property is of course (ii). If
a quasi-minimal space $X$ has a minimal subspace $Y$ then $Y$ embeds
into every subspace of~$X$. (Under these circumstances $X$ is said to
be $Y$-{\it saturated}.) For if $Z$ were a subspace into which $Y$ did
not embed then $Y$ would embed into every subspace of $Y$ and no
subspace of $Z$, so $Y$ and $Z$ would be totally incomparable. It
follows that $X$ has a minimal subspace if and only if $\cs$ has
cardinality one for every choice of $\cs$. If $\cs$ has cardinality
greater than one, in which case we shall say that $X$ is {\it strictly}
quasi-minimal, then $\cs$ is uncountable. To see this, suppose $\cs$
is countably infinite (the finite case is easier). We can pick a chain
$Y_1>Y_2>\dots$ of subspaces in $\cs$ such that every $Y\in\cs$
satisfies $Y>Y_n$ for some $n$. Now let $z_1<z_2<\dots$ with $z_n\in
Y_n$ and let $Z$ be the subspace spanned by $z_1,z_3,z_5,\dots\, $. Then
some $Y_n$ embeds into $Z$.  However, $Z$ is isomorphic to the
subspace spanned by
$z_{n+1},z_{n+2},\dots,z_{2n},z_{2n+1},z_{2n+3},z_{2n+5},\dots$ which
is a subspace of $Y_{n+1}$, contradicting the fact that $Y_n>Y_{n+1}$.
Finally, observe that a hereditarily indecomposable space is
quasi-minimal for trivial reasons. The definition will therefore only
interest us in the context of spaces with an unconditional basis.
 
We say that two block subspaces $Y$ and $Z$ of a space $X$ are {\it disjointly supported} if the support of every $y\in Y$ is disjoint
from the support of every\break $z\in Z$. We can now state and prove a
dichotomy for spaces with an unconditional basis. Later we shall
strengthen it a little.
\vglue-12pt
\proclaim{Theorem}    Let $X$ be a Banach space with an unconditional
basis. Then either $X$ has a quasi\/{\rm -}\/minimal subspace or $X$ has a 
subspace $Y$ such that no two disjointly supported subspaces of
$Y$ are isomorphic {\rm (}\/and hence any two disjointly supported subspaces
are totally incomparable\/{\rm ).}\/
\endproclaim
 
Although it is possible to prove Theorem 7.2 in one go using Theorem
4.1, it is perhaps of interest to show that the full strength of
Theorem 4.1 is by no means necessary. Instead we shall use Corollary
4.7 (that N-closed sets are weakly Ramsey) and a diagonalization.  In
order to do this, we must, as with Theorem 1.4, introduce some more
quantitative definitions. Given two block subspaces $Y,Z$ of a Banach
space $X$, we shall say that they are $C$-{\it block isomorphic} if
the linear map from $Y$ to $Z$ taking the normalized block basis of
$Y$ to the normalized block basis of $Z$ is a $C$-isomorphism.  We
shall say that a space $X$ is $C$-{\it quasi-minimal} if any two block
subspaces $Y,Z$ \pagebreak have further subspaces $V,W$ which are $C$-block
isomorphic. In the next lemma, when $A$ and $B$ are finite subsets
of $\N$, we write $A<B$ for the statement that $\max A<\min B$.
\vglue-12pt

\proclaim{Lemma}   Let $X$ be a Banach space with an unconditional
basis and let $C\ge 1$ and $\e>0$. Then either $X$ has a
$(C+\e)$\/{\rm -}\/quasi\/{\rm -}\/minimal subspace or $X$ has a subspace $V$ such that no two
disjointly supported block subspaces of $V$ are $C$\/{\rm -}\/block
isomorphic.
\endproclaim

{\it Proof}. For convenience we shall assume that $X$ has the property that
$\nm{x+y}>\nm x$ whenever $x<y$ are nonzero vectors in $X$. Standard
arguments tell us that every Banach space (with a given monotone
basis) is almost isometric to a space with this property. Suppose that
every block subspace $W$ of $X$ contains a pair of disjointly
supported $C$-block-isomorphic subspaces $Y,Z$. By passing to
subsequences we may assume that $Y$ and $Z$ are generated by
normalized block bases $\sq y$ and $\sq z$ such that $y_i<z_{i+1}$ and
$z_i<y_{i+1}$ for every $i\ge 1$. We know that for every $i\ge 1$ the
vectors $y_i$ and $z_i$ are disjointly supported linear combinations
of the vectors in the block basis of $W$. Therefore, we can find a
sequence $w_1<w_2<\dots$ of nonzero vectors in the unit ball of $W$
and a partition $A_1<A_2<\dots$ of $\N$ such that, if we let $y_i$ be
the sum of all $w_j$ such that $j\in A_i$ is odd, and $z_i$ be the sum
of all $w_j$ such that $j\in A_i$ is even, then the sequences
$y_1,y_2,\dots$ and $z_1,z_2,\dots$ are $C$-equivalent normalized
block bases. Let $\s$ be the set of sequences $\sq w$ for which such a 
partition of $\N$ exists. (Of course, strictly speaking we should
replace each $w_i$ by the pair $(w_i/\nm{w_i},\nm{w_i})$.)

Notice that, given our sequence $w_1,w_2,\dots$, the sets $A_i$, and
hence the vectors $y_i$ and $z_i$, are uniquely determined. For
example, $A_1=\{1,\dots,n\}$, where $n$ is the unique positive integer
such that the sum of all $w_j$ with $j\le n$ odd has norm 1 and so
does the sum of all $w_j$ with $j\le n$ even. Once $A_1$ has been
determined, the uniqueness of $A_2$ follows similarly, and so on. The
sequences $y_1,y_2,\dots$ and $z_1,z_2,\dots$ are $C$-equivalent if
and only if all their initial segments are $C$-equivalent. By the
uniqueness of the sets $A_i$, this is a property of initial segments
of the sequence $w_1,w_2,\dots\, $.  Moreover, if $w_1,w_2,\dots$ is a
sequence not in $\s$, either because we cannot find the sets $A_i$
making the vectors $y_i$ and $z_i$ normalized, or because for some $n$
the sequences $\seq y n$ and $\seq z n$ are not $C$-equivalent, then a
sufficiently small perturbation of $w_1,w_2,\dots$ will also not
be in $\s$. This shows that $\s$ is N-closed.

Since $\s$ is closed and large for $X$, Corollary 4.7 allows us, for
any\break $\D>0$, to pass to a subspace $V\subset X$ such that $\s_\D$ is
strategically large for $Y$. It is not hard to show that if $\D$ is
sufficiently small, then for any sequence $\sq w\in\s_\D$ there are
normalized sequences $\sq y$ and $\sq z$ which are 
$(C+\e)$-equivalent such that the $y_i$ are generated by
$w_1,w_3,\dots$ and the $z_i$ are generated by $w_2,w_4,\dots\, $.
Since $\s_\D$ is strategically large for $V$, if S chooses any
two subspaces $Y,Z$ of $V$ and plays the strategy $Y,Z,Y,Z,\dots$,
P can find a sequence $\sq w\in\s_\D$ with $w_i\in Y$ for $i$ odd
and $w_i\in Z$ for $i$ even. Choosing sequences $\sq y$ and $\sq z$
as above, we find block subspaces of $Y$ and $Z$ which are
$(C+\e)$-block isomorphic. \hfill\qed

\demo{Proof of Theorem {\rm 7.2}}  Suppose that $X$ has no
quasi-minimal subspace. Then obviously $X$ has no $m$-quasi-minimal
subspace for any positive integer $m$. Therefore by Lemma 7.3 we can
find a nested sequence of subspaces $W_1\supset W_2\supset\dots$ such
that no two disjointly supported block subspaces of $Y_m$ are
$(m-1)$-block-isomorphic. Let $w_1<w_2<\dots$ be a normalized block
basis with $w_m\in W_m$ for every $m$ and let $W$ be the subspace
generated by $w_1,w_2,\dots\, $. We claim that no two disjointly
supported subspaces of $W$ are isomorphic.

To see this, suppose that $Y$ and $Z$ are two disjointly supported
block subspaces of $W$. If they are isomorphic, then by standard
arguments there exist equivalent block bases $\sq y$ and $\sq z$
with every $y_i$ in $Y$ and every $z_i$ in $Z$ (but not necessarily 
generating the whole of $Y$ and $Z$). They must be $C$-equivalent
for some $C$, so let us choose a positive integer $m\ge C+2$, and
suppose (by removing the beginnings) that every $y_i$ and every
$z_i$ lies in $W_m$. We are not quite finished because the $y_i$
and $z_i$ are not necessarily normalized. However, the ratio
$\nm{y_i}/\nm{z_i}$ is bounded and bounded away from zero. 
Therefore, we can drop to subsequences such that this ratio
converges very quickly to a nonzero limit, and then a small
perturbation gives us normalized block bases which are
$(C+1)$-equivalent. This then contradicts our choice of $W_m$. 
\enddemo

It is clear that the dichotomy above is genuinely a dichotomy. If a
space even contains two disjointly supported totally incomparable
subspaces then it is not quasi-minimal. As we mentioned above, it is
possible to prove Theorem 7.2 more directly. If we let $\s_m$ be the
set $\s$ defined in the proof of Lemma 7.3 when $C=m$ and let
$\s=\bigcup_{m=1}^\infty$, then $\s$ is an ${\rm F}_\s$-set. If every
subspace contains two isomorphic disjointly supported further
subspaces, then the argument in the deduction of Theorem 7.2 from
Lemma 7.3 shows that $\s$ is large. Applying Theorem 4.1, we can find
a subspace $W$ where $\s_\D$ is strategically large which proves, as
above, that $W$ is quasi-minimal. 

Notice that the proof of Theorem 7.2 actually gave us something more:
every Banach space $X$ has a subspace $Y$ such that either $Y$ is
$C$-quasi-minimal for some $C$ or $Y$ does not contain two disjointly
supported isomorphic subspaces. In particular, every quasi-minimal
space has a subspace which is\break $C$-quasi-minimal for some $C$.

This last result is not unlike a consequence of our solution to
Banach's problem on homogeneous spaces. It follows from the fact that
a homogeneous space $X$ must be a Hilbert space that there exists a
constant $C$ such that all subspaces of $X$ are $C$-isomorphic.  (One
might call such a space\break $C$-homogeneous.) However, no direct proof of
this fact is known.

We now give some more information about the ``bad'' spaces in 
Theorem~7.2. First, we make a simple observation using a result of
Casazza [C]. (The proof can also be found in [G4].)

\proclaim{Lemma}   Let $X$ be a space with a basis such that{\rm ,} for any
block basis $y_1,z_1,y_2,z_2,\dots$ the block bases $y_1,y_2,\dots$
and $z_1,z_2,\dots$ are not equivalent. Then $X$ is isomorphic to no
proper subspace of itself {\rm (}\/and hence the same is true in any subspace
of $X${\rm ).}
\endproclaim

\proclaim{{C}orollary}  Let $X$ be a space with an unconditional basis
such that no two disjointly supported subspaces are isomorphic. Then
no subspace of $X$ is isomorphic to a further proper subspace.
\endproclaim

\demo{Proof} Clearly such a space satisfies the above criterion of Casazza.
\enddemo

\noindent We shall now strengthen Corollary 7.5 considerably. In the
proof below, we identify a set $A\subset\N$ with the obvious
projection on $X$ associated with $A$. The result basically states
that two block subspaces of a space $X$ satisfying the conditions of
Corollary 7.5 are only isomorphic if they are isomorphic for trivial
reasons.

\proclaim{Lemma}   Let $X$ be a space with an unconditional basis
such that no two disjointly supported subspaces are isomorphic{\rm ,} and
let $Y$ and $Z$ be subspaces of $X$ generated by equivalent block
bases. Let $T:Y\ra Z$ be the corresponding isomorphism. Then there is
an invertible diagonal operator $D:X\ra X$ such that $T-D|_Y$ and
$T^{-1}-D^{-1}|_Z$ are strictly singular.
\endproclaim

\demo{Proof} Let $\sq y$ be a normalized block basis generating $Y$ and let
$z_n=Ty_n$ for every $n$. Without loss of generality the basis of $X$
is 1-unconditional. Let $c,C$ be constants such that $c\nm
y\le\nm{Ty}\le C\nm y$ for every $y\in Y$. For each $n$ let $J_n$,
$K_n$ and $L_n$ be respectively the sets of integers $m$ such that
$0<2C|y_{nm}|\le|z_{nm}|$, $2C|y_{nm}|>|z_{nm}|>(c/2)|y_{nm}|>0$ and
$(c/2)|y_{nm}|\ge|z_{nm}|$. Let $J$, $K$ and $L$ be the unions of the
$J_n$, $K_n$ and $L_n$. Notice that these sets are disjoint.

Now let $D'$ be the unique diagonal operator on $KX$ sending $Ky_n$ to
$Kz_n$ for every $n$. Extend $D'$ to an invertible diagonal operator
$D$ on the whole of $X$ by mapping $x\in(1-K)X$ to $\l x$, where $\l$
is a constant between $c$ and $C$. We claim that $D$ has the property
required.

We show first that the map $U$ defined on $KY$ by $K y_n\mapsto
(J+L)y_n$ is continuous. If this were not so, then since $\sq {{Jy}}$
and $\sq {{Ly}}$ are not equivalent we would be able to find $y\in Y$
such that either $\nm{Jy}\ge 2\nm{(K+L)y}$ or
$\nm{Ly}>(4C/c)\nm{(J+K)y}$. Now since the basis of $X$ is
1-unconditional and $\sq z$ is a block basis, we have in the
first case
\begin{eqnarray*}
 \nm{Ty}&\ge&\nm{TJy}-\nm{T(K+L)y}>2C\nm{Jy}-C\nm{(K+L)y}\\
&\ge&
C(\nm{Jy}+\nm{(K+L)y})\ge C\nm y
\end{eqnarray*}
contradicting the fact that $\nm T\le C$. Similarly, in the
second case, we have
\begin{eqnarray*}
\nm{Ty}&\le&\nm{TLy}+\nm{T(J+K)y}\\
&\le&(c/2)\nm{Ly}+C\nm{(J+K)y}
<(3c/4)\nm{Ly}\le c\nm y
\end{eqnarray*}
contradicting the fact that $\nm{T^{-1}}\le c^{-1}$. Since $J$, $K$
and $L$ are disjoint, we also know that $U$ cannot be an isomorphism
on any subspace of $KY$, so it is strictly singular.

It follows that the sequences $\sq y$ and $\sq {{Ky}}$ are equivalent.
Thus the map $V$ defined on $KY$ by $Vy_n=(1-K)z_n$ is continuous
and therefore also strictly singular.
But $T=(D+V)K=D-D(J+L)+VK$ so $T-D$ is strictly singular as claimed.
The argument for $T^{-1}$ is similar.
\enddemo

  Putting all these results together, we get the following
theorem.   

\proclaim{Theorem} Let $X$ be an infinite\/{\rm -}\/dimensional
Banach space.  Then $X$ has a subspace $Y$ with one of the following
properties{\rm ,} which are mutually exclusive and all possible.
\begin{itemize}
\ritem{(1)} $Y$ is hereditarily indecomposable{\rm ,} and therefore {\rm (}\/by 
{\rm [F2])} every operator from a subspace $Z$ of $Y$ into $Y$ is a strictly
singular perturbation of a multiple of the inclusion map.
\ritem{(2)} $Y$ has an unconditional basis and every isomorphism
between block subspaces $W$ and $Z$ of $Y$ is a strictly singular
perturbation of the restriction of some invertible diagonal operator
on $Y$.
\ritem{(3)} $Y$ has an unconditional basis and is strictly 
quasi\/{\rm -}\/minimal.
\ritem{(4)} $Y$ has an unconditional basis and is minimal.
\end{itemize}

\endproclaim

How might Theorem 7.7 be extended? The obvious class to look at is
(3).  Notice that it is not at all clear that a strictly quasi-minimal
space need be isomorphic to a proper subspace. In fact, it is almost
certainly not even true. Indeed, a suggested counterexample, with a very
sketchy argument about why it was a counterexample, was given in [G3],
but so far nobody, the author included, has checked whether the
details can be filled in. If such a space exists, it is
likely to be the ``worst'' quasi-minimal space in the sense that the
only operators on the space and its subspaces are essentially those \pagebreak
guaranteed to exist by quasi-minimality. One can then ask what there is
to say about a space not containing one of these ``worst'' subspaces and
hope to divide (3) further.

The ultimate extension of Theorem 7.7 would be a sequence of possible
and mutually exclusive properties (1) to (k) such that, for any Banach
space $X$, there would exist $j$ and a subspace $Y\subset X$ such that
$Y$ was the ``worst'' sort of space satisfying property (j) in the
sense that the only operators on $Y$ were those guaranteed by property
(j), and such that property (k) was that of being isomorphic to some
$\ell_p$. Tsirelson's space ought to be an example of a typical space
having a property a little stronger than (3) but not as strong as
(4). Such a result would in a sense be a complete theory of the
structure of operators on subspaces of Banach spaces, as it would show
that every space had a subspace that was somehow ``generic'' of a
certain kind, and therefore could be described in some detail.

Unfortunately, it is not altogether obvious what the next step
is. Suppose $X$ is strictly quasi-minimal with an unconditional basis
and suppose that for every subspace $Y$ of $X$ there are isomorphic
subspaces $W$, $Z$ of $Y$ with $W$ properly contained in $Z$. Then no
subspace of $X$ satisfies the criterion of Casazza (Lemma 7.4) so in
some further subspace there is a winning strategy for P for producing
sequences $\sq x$ such that the sequence of odd-numbered vectors is
equivalent to the sequence of even-numbered ones. This, although a
stronger property than quasi-minimality, is rather artificial, and it
is not clear what of interest can be deduced from it. We can also find
a winning strategy in some subspace for sequences $\sq x$ such that
there is a partition of $\N$ into sets $A_1,A_2,\dots$ with
$\max(A_i)<\min(A_{i+1})$ with $y_i=\sum_{j\in A_i}x_j$ of norm 1 for
every $i$ and with $(x_i/\nm{x_i})_{i=1}^\infty$ and $\sq y$
equivalent block bases. Again, it is not clear how to use this fact.

It may well be that there is no tidy dichotomy connected with
isomorphisms to proper subspaces. If the following conjecture is true,
then it would certainly place a limit on any dichotomy that one might
wish to prove.
\vglue4pt
{\elevensc Conjecture 7.8}. {\it There exists a Banach space $X$ such
that every subspace $Y$ of $X$ has further subspaces $Z$ and $W$
such that $Z$ is isomorphic to no proper subspace of itself and
$W$ is isomorphic to its hyperplanes.}
\vglue6pt\advance\theoremcount by 1

 The difficulty in proving this conjecture is that the
known techniques for finding spaces that are not isomorphic to
any proper subspace all produce examples such that all subspaces
have the same property. A new idea is needed for ruling out 
operators on the whole space without doing the same for subspaces.

It might be a good idea to focus on the following rather general 
question. Recall the definitions connected with Lemma 7.1, that
$Y\le X$ if $Y$ is isomorphic to a subspace of $X$ and $Y\sim X$
if $Y\le X$ and $X\le Y$. 

\proclaim{Problem} Given a Banach space $X${\rm ,} let ${\cal P}(X)$
be the set of all equivalence classes of subspaces of $X${\rm ,} partially 
ordered as above. For which posets $P$ does there exist a Banach space
$X$ such that every subspace $Y$ of $X$ contains a further subspace $Z$
with ${\cal P}(Z)=P${\rm ?}
\endproclaim 

  Characterizing all such posets might be rather difficult, if
it involved constructing exotic Banach spaces. However, even a strong
necessary condition, perhaps proved using the methods of this paper,
would be interesting. 

Let us conclude this section with a vaguely stated problem.

\proclaim{Problem}\hskip-8pt Are there further nonartificial applications
of Theorem~{\rm 4.1?} In particular{\rm ,} is there any application which 
needs the full strength of the theorem\/{\rm ?}
\endproclaim

 Although I have searched unsuccessfully for such an
application, I still believe that Theorem 4.1 can be exploited
further. Many natural classes of sequences, such as the set of all
sequences that generate a subspace isomorphic to an $\ell_p$ space,
concern hereditary properties of subspaces. However, Theorem 4.1 has
nothing interesting to say about such properties, since it is trivial
that if $\s$ is a hereditary property of subspaces, then there is a
subspace such that either all its subspaces are in $\s$ or none of
them are. This indicates that any application of Theorem 4.1 is likely
to involve a ``clever'' choice of sequences, which gives me hope that
there are interesting applications that have not yet been found. As
for using the full strength of the theorem, this could mean applying
it to a set of sequences which is analytic but not Borel (see [Bos] for
some good examples of natural Banach-space properties which are
genuinely analytic), but it would also be good to see an application
which involved a strategy for S that was more complicated than simply
alternating between two subspaces.
 
\vglue-14pt
\section{Appendix}
\vglue-6pt
 
In this section we indicate briefly directions in which the results of
Sections 2 to 5 cannot be strengthened. The most obvious thing to
show is that not every set is weakly Ramsey. That follows from the
following statement, that there need not be a subspace in which the
game is even approximately determined. The proof, of course, uses the
axiom of choice, and in fact the continuum hypothesis as well,
although this is not so obviously necessary. (In fact, in [BLA1] it
is shown that Martin's axiom suffices.) It would be interesting
to know whether the axiom of determinacy implied that every set was
weakly Ramsey.

In order to apply the continuum hypothesis, it is essential that the
number of strategies for both players of the game is at most
$2^{\aleph_0}$.  However, this is not true with the definition of the
game given earlier, since both players have $2^{\aleph_0}$ possible
moves at each stage. Fortunately, this problem is not too
serious. First of all, one can restrict each of P's moves to a
countable dense set. (For example, insist that they are rational
linear combinations of the block basis presented by S.) To deal with
S, notice that it is not important for the game that S should declare
an entire block subspace all at once --- it is enough if S presents the
block basis one vector at a time, letting P choose a linear
combination of blocks at some finite stage. If in addition S is
restricted to a countable dense set, then the number of strategies for
both players is $2^{\aleph_0}$ as required. I am grateful to Joan
Bagaria and Jordi Lopez Abad for pointing out the necessity of this
remark.

\proclaim{Theorem}   Let $\D$ be the constant sequence
$(1/2,1/2,\dots)$. For every $X$ there exists a set $\s
\subset\Sigma(X)$ such that in no subspace $Y\subset X$ does $P$ have a
winning strategy for the game $\s_\D[Y]$ and in no subspace $Y\subset
X$ does $S$ have a winning strategy for the game $\s^{{\rm c}}[Y]$.
\endproclaim

\demo{Proof} Well order the set of all possible pairs consisting of a
pair $(Y,\phi)$, where $Y$ is a block subspace of $X$ and $\phi$ is
either a strategy for P or a strategy for S inside $Y$. Do this so
that each pair has countably many predecessors. We inductively define
for each ordinal $\a$ a sequence $s_\a$ as follows. Suppose the
$\a^{{\rm th}}$ pair in the well ordering is the pair $(Y_\a,\phi_\a)$
and that $\phi_\a$ is a strategy for P in $Y_\a$. Since there are only
countably many sequences $s_\b$ with $\b<\a$, S can easily find a
sequence of moves such that the sequence chosen by P is not within
$\D$ of any of these $s_\b$. In that case, let $s_\a$ be the sequence
resulting from such a play from S (with P playing the strategy
$\phi_\a$). Call such an $s_\a$ an S-{\it sequence}. Similarly, if
$\phi_\a$ is a strategy for S, then P can easily choose a sequence $s_\a$
against this strategy not within $\D$ of any $s_\b$ with $\b<\a$. Call
this sort of $s_\a$ a P-sequence. Note that, in particular, we have
guaranteed that no S-sequence is within $\D$ of any P-sequence.

Now let $\s$ be the set of all P-sequences. Given any strategy for
S and any subspace $Y$, P can play against this strategy and produce
a P-sequence. Thus, S has no strategy for any $\s^{{\rm c}}[Y]$.
Similarly, given any strategy for P and any subspace $Y$, S can
force P to produce an S-sequence, which is not in $\s_\D$. So
P does not have a winning strategy for $\s_\D[Y]$.
\enddemo

Recall the definition of an asymptotic set given earlier. A
fundamental fact in Banach space theory is that there exist
spaces $X$ such that the unit sphere $S(X)$ contains two
asymptotic sets $A$ and $B$ and $\d>0$ with $\nm{x-y}\ge\d$ for
every $x\in A$ and $y\in B$. It follows from work of Milman [Mi] 
that such sets exist in Tsirelson's space. A major breakthrough 
due to Odell and Schlumprecht [OS] was the discovery of such sets in
$\ell_p$ for $1\le p<\infty$. The results of [Mi], [OS] and [G5] in
fact show that a space has this property if and only if it does
not contain $c_0$ [OS, Corollary 1].

\proclaim{{C}orollary}  If $X$ does not contain $c_0${\rm ,} then there
exist subsets $\s$ and $\tau$ of $\Sigma(X)$, and $\G>0${\rm ,} such that $P$
has a winning strategy for every $\s[Y]$ and every $\tau[Y]$ but $P$ has
no winning strategy for any $(\s_{\G}\cap\tau_{\G})[Y]$ and $S$ has no
winning strategy for any $(\s\cap\tau)^{{\rm c}}[Y]$. 
\endproclaim

\demo{Proof} Let $\cal A$ and $\cal B$ be asymptotic sets such that
$\nm{x-y}\ge\d$ for every $x\in \cal A$ and $y\in \cal B$. Let $\rho$ be a set
given by Theorem 8.1. Let $\s$ be the set of all sequences $\sq x$
either in $\rho$ or such that $x_1\in \cal A$ and let $\tau$ be the same
with $\cal B$. Clearly P has winning strategies in every subspace $Y$ for
$\s[Y]$ and for $\tau[Y]$. Let $\G$ be a sequence with
$\gamma_1=\min(\d/2,1/3)$. Then
$\s_\G\cap\tau_\G=\rho_\G\subset\rho_\D$ and $\s\cap\tau=\rho$, so the
result follows from the property of $\rho$.  \enddemo

If $\s$ and $\tau$ are given by this corollary, then $\s_{\D/2}$ and
$\tau_{\D/2}$ are obviously completely weakly Ramsey, but their 
intersection is not even in an approximate sense completely weakly
Ramsey. 

\proclaim{Lemma} If $X$ does not contain $c_0$ then there exists a
large subset $\s\subset\Si(X)$ and $\G>0$ such that $\s_\G$ is
{\rm D-}\/$*$\/{\rm -}\/nowhere dense.
\endproclaim

\demo{Proof} Let $\cal A$ and $\G$ be as in the proof of the previous
corollary.  Let $\s$ be the set of sequences $\sq x$ with $x_n\in \cal
A$ for every $n$. Clearly $\s$ is large. However, given any basic open
set $[A;Z]$, pick $z\in Z\cap\cal B$. Then
$[(A,z);Z]\cap\s_\G\break =\emptyset$. \enddemo

\section{Recent developments}

As mentioned in the introduction, this paper is the final version of a
preprint [G3] that has been around for several years. I can offer no
excuse for this state of affairs, but am glad to remedy it now. I am
grateful to certain people for drawing my attention to mistakes in the
preprint. Valentin Ferenczi pointed out that the proof of Theorem 8 in
that paper (the statement that ${\rm D}$-open sets are weakly Ramsey) was
incorrect. (On the other hand, the proof in [G2] was correct.)  This
was later pointed out to me again by Joan Bagaria and Jordi Lopez
Abad, who also said that they found the proof of the main theorem
(corresponding to Theorem 4.1 in this paper) hard to follow.  In
response, I have made my presentation clearer in many respects. I have
also noticed and corrected a mistake in the proof of Lemma 9, which
corresponds to Lemma 4.12 here. (Lemma 4.11 does not have a
counterpart in the preprint, and was the missing step.) 

Recently Bagaria and Lopez Abad have given alternative proofs of the
main results of this paper, and have extended some of them [BLA1,2].
Their arguments are along similar lines to the ones here, but are
expressed in a language more familiar to logicians. For example, they
speak of dense generic filters when carrying out their
diagonalizations. Under additional set-theoretic hypotheses they
extend Theorem 4.1 to $\Sigma_2^1$-sets, that is, continuous images of
coanalytic sets.  Their work was carried out before I revised my
preprint.

\AuthorRefNames [BLA2]

\bye